# ON THE PERFORMANCE OF HIGH-ORDER FINITE ELEMENTS WITH RESPECT TO MAXIMUM PRINCIPLES AND THE NON-NEGATIVE CONSTRAINT FOR DIFFUSION-TYPE EQUATIONS

G. S. PAYETTE, K. B. NAKSHATRALA, AND J. N. REDDY

ABSTRACT. The main aim of this paper is to document the performance of $p$-refinement with respect to maximum principles and the non-negative constraint. The model problem is (steady-state) anisotropic diffusion with decay (which is a second-order elliptic partial differential equation). We considered the standard single-field formulation (which is based on the Galerkin formalism) and two least-squares-based mixed formulations. We have employed non-uniform Lagrange polynomials for altering the polynomial order in each element, and we have used $p = 1, \cdots, 10$. It will be shown that the violation of the non-negative constraint will not vanish with $p$-refinement for anisotropic diffusion. We shall illustrate the performance of $p$-refinement using several representative problems. The intended outcome of the paper is twofold. Firstly, this study will caution the users of high-order approximations about its performance with respect to maximum principles and the non-negative constraint. Secondly, this study will help researchers to develop new methodologies for enforcing maximum principles and the non-negative constraint under high-order approximations.

## 1. INTRODUCTION

Contaminant transport, chemical and biological remediation, and carbon-dioxide sequestration are some of the main engineering challenges of the 21st century. A large-scale implementation of any of these problems will have irreversible consequences on the environment, and can affect large geographical region and for long periods of time. In addressing these pressing issues, robust predictive numerical simulations have played and will continue to play an important role.

In these kind of problems, predicting the fate of chemical species is an important component, and diffusion is a dominant phenomenon. It should be noted that concentration (and certain other quantities such as absolute temperature, density, absolute pressure, saturation in multiphase systems) attains only non-negative values. A negative value for the concentration is unphysical. In a coupled reactive-transport numerical simulator, a negative value for the concentration of a species will result in an algorithmic failure. A robust numerical simulation should therefore preserve the physical and mathematical requirement of non-negativeness of species concentration.

1.1. **Maximum principles and high-order approximations.** Mathematically speaking, steady-state diffusion-type equations are elliptic partial differential equations, and are known to satisfy the so-called *maximum principles* [10]. The non-negative constraint can be obtained from maximum principles under certain assumptions on the input data. The discrete version of maximum principles are commonly referred to as discrete maximum principles (DMP). A study on discrete maximum principles investigates whether a given numerical formulation (that is, in a discrete setting) inherits the underlying maximum principles (which are satisfied in the continuous setting). A robust predictive numerical simulation should preserve various fundamental properties like the non-negativeness and maximum principles. Hence, one is particularly interested in the necessary and/or sufficient conditions under which a given numerical formulation satisfies discrete maximum principles and the non-negative constraint.

Many numerical formulations (under finite element, finite volume, and finite difference methodologies) have been developed. However, it should *not* be expected that any of these formulations will satisfy maximum principles and the non-negative constraint as there are no in-built mechanisms in these formulations to meet such constraints. Recently, researchers have proposed numerical methodologies for enforcing maximum principles and the non-negative constraint. For example, References [19, 20] have addressed the non-negative constraint under the finite volume method. Liska and Shashkov [21] have proposed a methodology to meet these constraints using the conservative finite difference technique. In References [23, 22], optimization-based techniques have been employed to meet these constraints under both mixed and single-field finite element formulations but have restricted their studies to low-order finite elements.

One of the early studies on discrete maximum principles is due to Varga [34], which was in the context of finite difference. An important work on discrete maximum principles with respect to the finite element method is by Ciarlet and Raviart [9]. In this paper, the authors have shown that acute-angle triangulation is a sufficient condition for a low-order approximation to satisfy maximum principles and the non-negative constraint. However, this condition is not sufficient under a high-order approximation. Other notable works on discrete maximum principles for high-order approximation are [14, 36, 31, 30, 35]. All these works considered one-dimensional problems. Herein we systematically study the performance of high-order approximations on several two-dimensional problems.

Also, earlier numerical works on discrete maximum principles concentrated on single-field formulations [22], and mixed formulations based on the variational multiscale formalism or lowest-order Raviart-Thomas spaces [23]. Herein, we shall also investigate the performance of mixed formulations based on the least-squares formalism.

1.2. **Least-squares formulations.** Least-squares variational methods often constitute an appealing alternative to the more popular weak formulation based on the Galerkin formalism in developing efficient and robust finite element models. In a least-squares-based formulation a non-physical *least-squares* functional is defined in terms of the sum of the squares of appropriate norms of the governing partial differential equation residuals. Since the least-squares functional is by definition positive and convex and it naturally follows that the minimizer of the least-squares functional coincides with the exact solution of the original set of partial differential equations. The finite element model associated with the least-squares formulation is constructed via direct minimization of the least-squares functional with respect to the trial space associated with the finite element discretization. Some representative references on least-squares formulations are [4, 16].

1.3. **Main contributions of this paper.** The main aim of this paper is to document the performance of $p$-refinement for solving diffusion-type equations with respect to maximum principles and the non-negative constraint. We illustrate the performance on various computational grids and using different canonical problems. We consider three weak formulations: the standard single-field formulation and two least-squares-based formulations. We illustrate the extent to which these formulations violate maximum principles and the non-negative constraint under high-order approximations. We do not, however, provide a methodology for enforcing maximum principles and the non-negative constraint under high-order approximations. We also show that the performance of least-squares formulations with respect to the non-negative constraint depends on the choice of the weight in the inner product. The present work has two main purposes. First, users of high-approximations will be aware of their performance with respect to maximum principles and the non-negative constraint. Second, it will help researchers to develop methodologies for enforcing maximum principles and the non-negative under high-order approximations.

1.4. **An outline of the paper.** The remainder of this paper is organized as follows. In Section 2, we present the governing equations of anisotropic diffusion with decay. We shall also discuss the classical maximum principle, and its consequences (in particular, the non-negative constraint). In Section 3, we present the weak formulations, and in Section 4 we discuss high-order spectral approximations. In Section 5, we shall show the performance of these numerical formulations using several canonical problems, and conclusions are drawn in Section 6.

We shall employ the following symbolic notation in this paper. To distinguish vectors in the continuum setting from vectors in the finite element context, we shall employ lower case boldface normal letters for the former, and lower case boldface italic letters for the later. For example, **x** is used to denote a spatial position vector, and $\boldsymbol{c}$ is used to denote a finite element vector containing nodal concentrations. To distinguish a second-order tensor from a matrix, we shall denote second-order continuum tensors using upper case boldface normal letters, and shall denote matrices using upper case boldface italic letters. For example, $\mathbf{D}$ is used to denote the diffusivity tensor, and $\boldsymbol{K}$ is used to denote stiffness matrix. Throughout this paper, repeated indices do not imply summation. (That is, Einstein's summation convention is *not* employed.) We shall denote the set of natural numbers as $\mathbb{N}$, and the set of real numbers as $\mathbb{R}$. Other notational conventions are introduced as needed.

## 2. GOVERNING EQUATIONS: DIFFUSION WITH DECAY

In this paper, we assume that $\Omega$ is an open bounded subset of $\mathbb{R}^{nd}$, where "*nd*" denotes the number of spatial dimensions. The boundary is denoted by $\partial \Omega := \bar{\Omega} - \Omega$, where a superposed bar denotes the set closure. The boundary is divided into two parts: $\Gamma^{\mathrm{D}}$ and $\Gamma^{\mathrm{N}}$, where $\Gamma^{\mathrm{D}}$ is that part of the boundary on which Dirichlet boundary conditions are prescribed, and $\Gamma^{\mathrm{N}}$ is the part of the boundary on which Neumann boundary conditions are prescribed. For well-posedness, we have $\Gamma^{\mathrm{D}} \cup \Gamma^{\mathrm{N}} = \partial \Omega$ and $\Gamma^{\mathrm{D}} \cap \Gamma^{\mathrm{N}} = \varnothing$. A spatial point is denoted by $\mathbf{x} \in \bar{\Omega}$. The gradient and divergence operators with respect to $\mathbf{x}$ are, respectively, denoted by $\mathrm{grad}[\cdot]$ and $\mathrm{div}[\cdot]$. Herein, we shall consider the anisotropic diffusion of a chemical species, and allow the decay of the chemical species. Let us denote the concentration of the chemical species by $c(\mathbf{x})$. We present the governing equations in the standard divergence form, which take the following form:

$$\alpha(\mathbf{x})c(\mathbf{x}) - \mathrm{div}[\mathbf{D}(\mathbf{x})\mathrm{grad}[c]] = f(\mathbf{x}) \quad \text{in } \Omega \tag{1a}$$

$$c(\mathbf{x}) = c^{\mathrm{p}}(\mathbf{x}) \quad \text{on } \Gamma^{\mathrm{D}} \tag{1b}$$

$$\mathbf{n}(\mathbf{x}) \cdot \mathbf{D}(\mathbf{x})\mathrm{grad}[c] = t^{\mathrm{p}}(\mathbf{x}) \quad \text{on } \Gamma^{\mathrm{N}} \tag{1c}$$

where $\alpha(\mathbf{x}) \geq 0$ is the decay coefficient, $\mathbf{D}(\mathbf{x})$ denotes the diffusivity tensor, $f(\mathbf{x})$ is the volumetric source/sink, $c^{\mathrm{p}}(\mathbf{x})$ is the prescribed concentration, $t^{\mathrm{p}}(\mathbf{x})$ is the prescribed flux, and $\mathbf{n}(\mathbf{x})$ is the unit outward normal vector on the boundary. The mathematical model given by equations (1a)–(1c) frequently arises in Mathematical Physics (see the discussion in [22, Introduction]).

We shall assume that the decay coefficient is bounded above, which means that there exists a real constant $\alpha_0 < +\infty$ such that we have

$$\alpha(\mathbf{x}) < \alpha_0 \quad \forall \mathbf{x} \in \Omega \tag{2}$$

(Note that we have already assumed that the decay coefficient is non-negative.) The diffusivity tensor is assumed to be symmetric, bounded above and uniformly elliptic. That is, there exists two constants $0 < \lambda_1 \leq \lambda_2 < +\infty$ such that

$$\lambda_1 \mathbf{y}^T \mathbf{y} \leq \mathbf{y}^T \mathbf{D}(\mathbf{x}) \mathbf{y} \leq \lambda_2 \mathbf{y}^T \mathbf{y} \quad \forall \mathbf{x} \in \Omega, \forall \mathbf{y} \in \mathbb{R}^{nd} \tag{3}$$

It should be noted that uniform ellipticity condition is a stronger requirement than demanding that the matrix corresponding to tensor $\mathbf{D}(\mathbf{x})$ being positive definite at every spatial point in the domain. The system of equations (1a)–(1c) is a second-order elliptic partial differential equation with Dirichlet and Neumann boundary conditions. It is well-known that the above partial differential equation satisfies the so-called *maximum principles* under certain regularity assumptions on the domain and on the input data ($\mathbf{D}(\mathbf{x})$, $\alpha(\mathbf{x})$, $f(\mathbf{x})$, $c^{\mathrm{p}}(\mathbf{x})$, and $t^{\mathrm{p}}(\mathbf{x})$).

2.1. **Maximum principles and the non-negative constraint.** There are different kinds of maximum principles for second-order elliptic partial differential equations available in the literature (for example, see the discussion in Han and Lin [12]). Herein, we consider two maximum principles by E. Hopf [15]. (Also see the commentary on Hopf's paper by J. Serrin [29].) The first maximum principle is for pure diffusion without decay (that is, $\alpha(\mathbf{x}) = 0$), and the second maximum principle allows the possibility of decay (that is, $\alpha(\mathbf{x}) \geq 0$). For $k \in \mathbb{N}$, we shall use $C^k(\Omega)$ to denote the set of functions having all derivatives up to the order $k$ continuous on $\Omega$, and $C^0(\bar{\Omega})$ to denote the space of all continuous functions on $\Omega$ that can be continuously extended to the boundary.

**Theorem 1** (A maximum principle for diffusion equation). *Let $c(\mathbf{x}) \in C^2(\Omega) \cap C^0(\bar{\Omega})$ satisfy the following differential inequality*

$$-\mathrm{div}[\mathbf{D}(\mathbf{x})\mathrm{grad}[c]] = f(\mathbf{x}) \leq 0$$

*with $\mathbf{D}(\mathbf{x})$ is uniformly elliptic, bounded above and continuously differentiable. Then, we have*

$$\max_{\mathbf{x} \in \Omega} c(\mathbf{x}) \leq \max_{\mathbf{x} \in \bar{\Omega}} c(\mathbf{x})$$

**Theorem 2** (A maximum principle for diffusion with decay). *Let $c(\mathbf{x}) \in C^2(\Omega) \cap C^0(\bar{\Omega})$ satisfy the following differential inequality*

$$\alpha(\mathbf{x})c(\mathbf{x}) - \mathrm{div}[\mathbf{D}(\mathbf{x})\mathrm{grad}[c]] = f(\mathbf{x}) \leq 0$$

*with $\alpha(\mathbf{x}) \geq 0$; $\mathbf{D}(\mathbf{x})$ is uniformly elliptic, bounded above and continuously differentiable; and $\alpha(\mathbf{x})$ is bounded. Then, we have*

$$\max_{\mathbf{x} \in \Omega} c(\mathbf{x}) \leq \max_{\mathbf{x} \in \bar{\Omega}} \{c(\mathbf{x}), 0\}$$

Mathematical proofs to the above two theorems can be found in Gilbarg and Trudinger [10]. In both theorems, we have assumed volumetric sink (i.e., $f(\mathbf{x}) \leq 0$). If we have volumetric source (i.e., $f(\mathbf{x}) \geq 0$) then the "max" has to be replaced by "min." Put differently, for volumetric source, the minimum occurs on the boundary in the case of pure diffusion, and the non-negative occurs on the boundary in the case of diffusion with decay. Similarly, if $f(\mathbf{x}) = 0$, then both the maximum and minimum occur on the boundary in the case of pure diffusion, and the non-negative maximum and non-negative minimum occur on the boundary in the case of diffusion with decay. If there are regions of both source and sink in the domain then the above theorems do not directly apply. The discrete version of maximum principles is commonly referred to as discrete maximum principles.

There are several studies on discrete maximum principles, and some representative ones are [8, 18, 13, 19, 21, 23, 22]. Some of these studies (particularly the ones that derive necessary and sufficient conditions to meet the non-negative constraint and maximum principles, for example, references [18]) hinge on the following mathematical arguments. After spatial discretization using the finite element method (or, for that matter, using the finite volume method or finite difference method) one obtains a system of linear equations of the following form:

$$\boldsymbol{Ax} = \boldsymbol{b} \tag{4}$$

where $\boldsymbol{A}$ is called the coefficient matrix, and $\boldsymbol{b}$ is called the force matrix. Under the low-order approximation (i.e., $p = 1$), it can be shown that $\boldsymbol{b} \succeq \boldsymbol{0}$ if the forcing function $f(\mathbf{x}) \geq 0$. (Herein, we have used the symbol $\succeq$ to denote component-wise inequality. That is, $\boldsymbol{b} \succeq \boldsymbol{0}$ means that $b_i \geq 0$.) If $\boldsymbol{b} \succeq \boldsymbol{0}$, a sufficient condition for $\boldsymbol{x} \succeq \boldsymbol{0}$ is requiring the matrix $\boldsymbol{A}$ to be monotone (which means that the inverse of $\boldsymbol{A}$ has all non-negative entries). Usually, a stronger condition requiring that the matrix $\boldsymbol{A}$ is an M-matrix is commonly employed [18]. (There are different ways of defining an M-matrix, and for further details see references [2, 28].) Restrictions are then placed (typically on the mesh) to make $\boldsymbol{A}$ to be an M-matrix. However, as discussed in references [23] and [22], such a procedure will work only for isotropic diffusion, and is not sufficient for anisotropic diffusion. Several examples are presented in these two references to illustrate this conclusion.

For high-order approximations, another complexity arises due to the fact that even if $f(\mathbf{x}) \geq 0 \ \forall \mathbf{x} \in \Omega$, its $L_2$ projection onto the polynomial space need not be non-negative. That is, even if $f(\mathbf{x}) \geq 0$, the condition $\boldsymbol{b} \succeq \boldsymbol{0}$ need not hold under high-order approximations. Hence, the whole analysis requiring that $\boldsymbol{A}$ to be a monotone or an M-matrix will not guarantee that the non-negative constraint and maximum principles under high-order approximations will be met. In this paper, we carefully assess the performance of high-order approximations with respect to maximum principles and the non-negative constraint.

**Remark 3.** *Some test problems in Section 5 do not have classical solutions in the sense that $c(\mathbf{x}) \in C^2(\Omega) \cap C^0(\bar{\Omega})$, and hence, theorems 1 and 2 do not directly apply. In the literature, however, one can find maximum principles for weak solutions posed under weaker regularity assumptions. For example, see references [25, 7, 5]. Such a detailed discussion is beyond the scope of this paper, and is not central to our presentation.*

**Remark 4.** *Some prior numerical works have addressed the non-negative constraint and maximum principles for models similar to equation (1) on general computational grids. Nakshatrala and*

Valocchi [23] have considered discrete maximum principles for pure diffusion and two mixed formulations. Specifically, they considered the variational multiscale formulation and the lowest-order Raviart-Thomas element. They provided a methodology for enforcing maximum principles and the non-negative constraint that works only for low-order finite elements. Nagarajan and Nakshatrala [22] have considered discrete maximum principles for the diffusion equation with decay under the standard single-field formulation. They also restricted their studies to low-order finite element approximations. Liska and Shashkov [21] have addressed discrete maximum principles for pure diffusion using conservative finite difference techniques. The present study differs from previous works in two ways: use and performance study of (1) high-order approximations, and (2) least-squares finite element models.

## 3. WEAK FORMULATIONS

In this paper we shall consider three weak formulations. The first is the classical single-field formulation, which is based on the Galerkin formalism. The second and third are mixed least-squares-based formulations. For completeness, we briefly outline these formulations. To this end, let us introduce the following function spaces:

$$\mathcal{U} := \left\{ c(\mathbf{x}) \in H^1(\Omega) \mid c(\mathbf{x}) = c^{\mathrm{p}}(\mathbf{x}) \text{ on } \Gamma^{\mathrm{D}} \right\} \tag{5a}$$

$$\mathcal{W} := \left\{ w(\mathbf{x}) \in H^1(\Omega) \mid w(\mathbf{x}) = 0 \text{ on } \Gamma^{\mathrm{D}} \right\} \tag{5b}$$

$$\mathcal{Q} := \left\{ \mathbf{q}(\mathbf{x}) \mid \mathbf{q}(\mathbf{x}) \in (H^1(\Omega))^{nd} \right\} \tag{5c}$$

where $H^1(\Omega)$ is a standard Sobolev space [32]. (Recall that "$nd$" denotes the number of spatial dimensions.)

### 3.1. The classical single-field formulation.
The weak formulation based on the Galerkin formalism for the governing equations (1a)-(1c) can be stated as follows: Find $c(\mathbf{x}) \in \mathcal{U}$ such that we have

$$\mathcal{B}(w(\mathbf{x}); c(\mathbf{x})) = \mathcal{F}(w(\mathbf{x})) \quad \forall w(\mathbf{x}) \in \mathcal{W} \tag{6}$$

where the bilinear form $\mathcal{B}(w(\mathbf{x}); c(\mathbf{x}))$ and the linear functional $\mathcal{F}(w(\mathbf{x}))$ are, respectively, defined as follows:

$$\mathcal{B}(w(\mathbf{x}); c(\mathbf{x})) := \int_\Omega \alpha(\mathbf{x}) w(\mathbf{x}) c(\mathbf{x}) \, \mathrm{d}\Omega + \int_\Omega \mathrm{grad}[w] \cdot \mathbf{D}(\mathbf{x}) \mathrm{grad}[c] \, \mathrm{d}\Omega \tag{7a}$$

$$\mathcal{F}(w(\mathbf{x})) := \int_\Omega w(\mathbf{x}) f(\mathbf{x}) \, \mathrm{d}\Omega + \int_{\Gamma^{\mathrm{N}}} w(\mathbf{x}) t^{\mathrm{p}}(\mathbf{x}) \, \mathrm{d}\Gamma \tag{7b}$$

One can easily verify that the above weak formulation (23) is equivalent to minimizing the following functional:

$$\mathcal{J}(c(\mathbf{x})) := \frac{1}{2} \int_\Omega \Big(\alpha(\mathbf{x})c^2(\mathbf{x}) + \mathrm{grad}[c] \cdot \mathbf{D}(\mathbf{x})\mathrm{grad}[c]\Big) \mathrm{d}\Omega - \int_\Omega c(\mathbf{x})f(\mathbf{x})\,\mathrm{d}\Omega$$
$$- \int_{\Gamma^\mathrm{N}} c(\mathbf{x})t^\mathrm{p}(\mathbf{x})\,\mathrm{d}\Gamma = \frac{1}{2}\mathcal{B}(c(\mathbf{x});c(\mathbf{x})) - \mathcal{F}(c(\mathbf{x})) \qquad (8)$$

Hence, the finite element approximation inherits the best approximation property with respect to the norm $\|c(\mathbf{x})\|_E = \sqrt{\frac{1}{2}\mathcal{B}(c(\mathbf{x});c(\mathbf{x}))}$.

3.2. **Least-squares-based formulations.** Least-squares variational procedures possess many attractive mathematical as well as practical computational attributes. In particular, the least-squares method always invokes a minimization principle whose minimizer coincides with the solution of the governing partial differential equations. As a result, the discrete numerical solution always possesses the best approximation property with respect to a well-defined norm (i.e., the energy norm of the functional). When this norm is equivalent to a standard norm (such as a norm from an appropriate Sobolev space) ideal convergence rates of the finite element solution may be established. That the least-squares formulation is always based on a minimization principle insures a robust setting that is often lacking in the associated weak form Galerkin model. In addition, the bilinear form resulting from invocation of the minimization principle in least-squares formulations will always be symmetric and positive definite, and restrictive compatibility conditions such as the *inf-sup* condition will never arise. As a result, providing a conforming discretization and well-posed boundary conditions, the discrete variational formulation will always possesses a unique solution. Furthermore, the symmetry and positive-definiteness of the resulting coefficient matrix may be exploited directly in the solution process of the linear system of equations. In particular, direct solvers may employ a sparse Cholesky decomposition of the coefficient matrix, while iterative solvers may utilize the robust preconditioned conjugate gradient algorithm. Due to symmetry, only half of the coefficient matrix need to be actually calculated and stored in memory.

Least-squares formulations are not without their own deficiencies and in many cases there are fixes. For example, unlike weak form Galerkin formulations where regularity requirements of the finite element spaces are weakened (by invoking Green's identities), least-squares formulations require higher regularity of the finite element spaces (dictated by the order of the governing partial differential equations). Higher regularity requirements negatively affect the condition number of the coefficient matrix and also the continuity requirement of the solution across element boundaries. High regularity requirements may be avoided by constructing the least-squares finite element model in terms of an equivalent lower-order system by the introduction of additional independent variables. Such a mixed formulation allows for the use of standard Lagrange interpolation functions but also produces a large set of global equations to be solved. Such a formulation is often quite valuable, however, since the auxiliary variables are typically physical quantities of interest such as the flux.

A least-squares-based formulation will be based on the minimization of a least-squares functional, which is constructed from the sum of the squares of appropriate norms of the residuals of the governing partial differential equations and the Neumann boundary condition. Although it is certainly possible to construct a least-squares finite element model based on equations (1a)-(1c) directly, such an approach requires a higher degree of regularity in the finite element solution such as $c(\mathbf{x}) \in H^2(\Omega)$, where $H^2(\Omega)$ is another standard Sobolev space on $\Omega$ [32]. Herein, we recast the governing equations into the following equivalent first-order system:

$$\alpha(\mathbf{x})c(\mathbf{x}) + \text{div}[\mathbf{q}] = f(\mathbf{x}) \quad \text{in } \Omega \tag{9a}$$

$$\mathbf{q}(\mathbf{x}) = -\mathbf{D}(\mathbf{x})\text{grad}[c] \quad \text{in } \Omega \tag{9b}$$

$$c(\mathbf{x}) = c^{\text{p}}(\mathbf{x}) \quad \text{on } \Gamma^{\text{D}} \tag{9c}$$

$$-\mathbf{n}(\mathbf{x}) \cdot \mathbf{q}(\mathbf{x}) = t^{\text{p}}(\mathbf{x}) \quad \text{on } \Gamma^{\text{N}} \tag{9d}$$

where $\mathbf{q}(\mathbf{x})$ is the diffusive flux vector field. We shall treat both $c(\mathbf{x})$ and $\mathbf{q}(\mathbf{x})$ as independent variables, and employ the standard Lagrange finite element interpolation functions for both of these field variables. The price incurred will be an increase in the overall size of the system of equations that one needs to solve to obtain a numerical solution. The additional expense is not unwarranted as the flux is an important quantity of interest in many engineering applications. As a result, one obtains the flux directly in the solution process as opposed to the usual way of obtaining it at the post-processing stage of a finite element analysis.

**Remark 5.** *For stability reasons, some least-squares formulations for diffusion-type equations augment the first order form of the governing partial differential equations (i.e., equations (9a)–(9d)) with the following seemingly redundant expression:*

$$\text{curl}[\mathbf{q}] = \mathbf{0} \tag{10}$$

*For example, see Reference [4]. It is important to note that $\mathbf{q}$ satisfies the above expression, in general, only if the medium is homogeneous (that is, the diffusivity tensor is independent of $\mathbf{x}$) and isotropic. In this work, we place no such restrictions on $\mathbf{D}(\mathbf{x})$, and hence do not utilize the above expression in construction of the least-squares functional.*

We shall define the norm $\|\cdot\|_K$ for scalar and vector functions defined over $K$ as follows:

$$\|u(\mathbf{x})\|_K^2 := \int_K u^2(\mathbf{x})\, \mathrm{d}K \tag{11a}$$

$$\|\mathbf{u}(\mathbf{x})\|_K^2 := \int_K \mathbf{u}(\mathbf{x}) \cdot \mathbf{u}(\mathbf{x})\, \mathrm{d}K \tag{11b}$$

In this paper we investigate the performance of two mixed least-squares-based formulations, and denote them as "LS1 formulation" and "LS2 formulation". For both these formulations, the least-squares functional can be constructed as follows:

$$\begin{aligned}
\mathcal{J}_{\mathrm{LS}}\left(c(\mathbf{x}), \mathbf{q}(\mathbf{x})\right) :=\ & \frac{1}{2} \|\beta(\mathbf{x})\left(\alpha(\mathbf{x})c(\mathbf{x}) + \mathrm{div}\,[\mathbf{q}] - f(\mathbf{x})\right)\|_\Omega^2 \\
& + \frac{1}{2} \|\mathbf{A}(\mathbf{x})\left(\mathbf{q}(\mathbf{x}) + \mathbf{D}(\mathbf{x})\mathrm{grad}[c]\right)\|_\Omega^2 \\
& + \frac{1}{2} \|\mathbf{q}(\mathbf{x}) \cdot \mathbf{n}(\mathbf{x}) + t^{\mathrm{p}}(\mathbf{x})\|_{\Gamma^{\mathrm{N}}}^2
\end{aligned} \tag{12}$$

where the second-order tensor $\mathbf{A}(\mathbf{x})$ and the scalar function $\beta(\mathbf{x})$ are, respectively, defined as follows:

$$\mathbf{A}(\mathbf{x}) = \begin{cases} \mathbf{I} & \text{LS1 formulation} \\ \mathbf{D}^{-1/2}(\mathbf{x}) & \text{LS2 formulation} \end{cases} \tag{13a}$$

$$\beta(\mathbf{x}) = \begin{cases} 1 & \text{LS1 formulation} \\ \left.\begin{array}{ll} 1 & \text{if } \alpha(\mathbf{x}) = 0 \\ \alpha^{-1/2}(\mathbf{x}) & \text{if } \alpha(\mathbf{x}) \neq 0 \end{array}\right\} & \text{LS2 formulation} \end{cases} \tag{13b}$$

where $\mathbf{I}$ is the second-order identity tensor. Recall that $\alpha(\mathbf{x}) \geq 0$, and hence $\beta(\mathbf{x})$ is well-defined. Also, note that $\mathbf{D}(\mathbf{x})$ is a positive definite tensor, and by square root theorem $\mathbf{D}^{-1/2}(\mathbf{x})$ is also well-defined [11]. The variational principle associated with both the least-squares formulations can be stated as follows: Find $c(\mathbf{x}) \in \mathcal{U}$ and $\mathbf{q}(\mathbf{x}) \in \mathcal{Q}$ such that we have

$$\mathcal{J}_{\mathrm{LS}}\left(c(\mathbf{x}), \mathbf{q}(\mathbf{x})\right) \leqslant \mathcal{J}_{\mathrm{LS}}\left(\tilde{c}(\mathbf{x}), \tilde{\mathbf{q}}(\mathbf{x})\right) \quad \forall \tilde{c}(\mathbf{x}) \in \mathcal{U},\ \tilde{\mathbf{q}}(\mathbf{x}) \in \mathcal{Q} \tag{14}$$

The first-order optimality condition demands that the first variation of $\mathcal{J}_{\mathrm{LS}}\left(c(\mathbf{x}), \mathbf{q}(\mathbf{x})\right)$ be identically zero. The corresponding weak statement can be written as follows: Find $c(\mathbf{x}) \in \mathcal{U}$ and $\mathbf{q}(\mathbf{x}) \in \mathcal{Q}$ such that we have

$$\mathcal{B}_{\mathrm{LS}}\left(w(\mathbf{x}), \mathbf{w}(\mathbf{x}); c(\mathbf{x}), \mathbf{q}(\mathbf{x})\right) = \mathcal{F}_{\mathrm{LS}}\left(w(\mathbf{x}), \mathbf{w}(\mathbf{x})\right) \quad \forall w(\mathbf{x}) \in \mathcal{W},\ \mathbf{w}(\mathbf{x}) \in \mathcal{Q} \tag{15}$$

where the bilinear form $\mathcal{B}_{\text{LS}}(w(\mathbf{x}), \mathbf{w}(\mathbf{x}); c(\mathbf{x}), \mathbf{q}(\mathbf{x}))$ and the linear functional $\mathcal{F}_{\text{LS}}(w(\mathbf{x}), \mathbf{w}(\mathbf{x}))$ are, respectively, defined as follows:

$$\mathcal{B}_{\text{LS}}(w(\mathbf{x}), \mathbf{w}(\mathbf{x}); c(\mathbf{x}), \mathbf{q}(\mathbf{x})) := \int_{\Omega} \beta^2(\mathbf{x}) \left(\alpha(\mathbf{x})w(\mathbf{x}) + \text{div}[\mathbf{w}]\right) \left(\alpha(\mathbf{x})c(\mathbf{x}) + \text{div}[\mathbf{q}]\right) \, d\Omega$$

$$+ \int_{\Omega} (\mathbf{w}(\mathbf{x}) + \mathbf{D}(\mathbf{x})\text{grad}[w]) \cdot \mathbf{A}^2(\mathbf{x}) (\mathbf{q}(\mathbf{x}) + \mathbf{D}(\mathbf{x})\text{grad}[c]) \, d\Omega$$

$$+ \int_{\Gamma^{\text{N}}} \mathbf{w}(\mathbf{x}) \cdot \mathbf{n}(\mathbf{x}) \, \mathbf{q}(\mathbf{x}) \cdot \mathbf{n}(\mathbf{x}) \, d\Gamma \tag{16a}$$

$$\mathcal{F}_{\text{LS}}(w(\mathbf{x}), \mathbf{w}(\mathbf{x})) := \int_{\Omega} \beta^2(\mathbf{x}) \left(\alpha(\mathbf{x})w(\mathbf{x}) + \text{div}[\mathbf{w}]\right) f(\mathbf{x}) \, d\Omega - \int_{\Gamma^{\text{N}}} \mathbf{w}(\mathbf{x}) \cdot \mathbf{n}(\mathbf{x}) \, t^{\text{P}}(\mathbf{x}) \, d\Gamma \tag{16b}$$

## 4. HIGH-ORDER APPROXIMATIONS

As shown in the previous section, a variational formulation (either Galerkin or least-squares formulations) of a general boundary value problem may be stated as follows: Find $\mathbf{u}(\mathbf{x}) \in \mathcal{V}$ such that we have

$$\mathcal{B}(\mathbf{w}(\mathbf{x}); \mathbf{u}(\mathbf{x})) = \mathcal{F}(\mathbf{w}(\mathbf{x})) \quad \forall \mathbf{w}(\mathbf{x}) \in \tilde{\mathcal{V}} \tag{17}$$

where $\mathcal{B}(\mathbf{w}(\mathbf{x}); \mathbf{u}(\mathbf{x}))$ is a bilinear form, $\mathcal{F}(\mathbf{w}(\mathbf{x}))$ is a linear form, and $\mathcal{V}$ and $\tilde{\mathcal{V}}$ are appropriate function spaces. The quantity $\mathbf{u}(\mathbf{x})$ represents the set of independent variables (associated with the variational boundary value problem), and $\mathbf{w}(\mathbf{x})$ represents the corresponding weighting function.

In the finite element method we restrict the solution space to a finite dimensional sub-space $\mathcal{V}^{hp}$ of the infinite dimensional space $\mathcal{V}$, and the weighting functions to a finite dimensional sub-space $\tilde{\mathcal{V}}^{hp} \subset \tilde{\mathcal{V}}$. As a result, in the discrete case we seek to find $\mathbf{u}_{hp}(\mathbf{x}) \in \mathcal{V}^{hp}$ such that we have

$$\mathcal{B}(\mathbf{w}_{hp}(\mathbf{x}); \mathbf{u}_{hp}(\mathbf{x})) = \mathcal{F}(\mathbf{w}_{hp}(\mathbf{x})) \quad \forall \mathbf{w}_{hp}(\mathbf{x}) \in \tilde{\mathcal{V}}^{hp} \tag{18}$$

We assume that the domain $\bar{\Omega} \subset \mathbb{R}^{nd}$ is discretized into a set of NE non-overlapping sub-domains $\bar{\Omega}^e$, called finite elements, such that $\bar{\Omega} \approx \bar{\Omega}^{hp} = \bigcup_{e=1}^{\text{NE}} \bar{\Omega}^e$. The geometry of each $\bar{\Omega}^e$ is characterized using the standard isoparametric bijective mapping from $\bar{\Omega}^e$ to the master element $\hat{\Omega}^e$. In the present study we restrict the classes of elements considered to lines in $\mathbb{R}^1$, four side quadrilateral elements in $\mathbb{R}^2$ and six face brick elements in $\mathbb{R}^3$ (although numerical results are presented for $nd = 1$ and 2 only). As a result we can simply define the master element as $\hat{\Omega}^e = [-1, +1]^{nd}$. The natural coordinates associated with $\hat{\Omega}^e$ (when $nd = 3$) are defined as $\boldsymbol{\xi} = (\xi, \eta, \zeta)$. In this work we utilize a family of finite elements constructed using high polynomial order interpolation functions. The quantities $h$ and $p$ appearing in the definition of the sub-space $\mathcal{V}^{hp}$ imply that the discrete solution may be refined by either increasing the number of elements in $\bar{\Omega}^{hp}$ ($h$-refinement), increasing the polynomial order of the approximate solution within each element $\bar{\Omega}^e$ ($p$-refinement) or through an appropriate and systematic combination of both $h$-refinement and $p$-refinement.

Within a typical finite element a given variable, such as the species concentration $c(\mathbf{x})$, may be approximated by the following formula

$$c(\mathbf{x}) \approx c_{hp}(\mathbf{x}) = \sum_{i=1}^{n} c_i^e \psi_i(\boldsymbol{\xi}) \qquad \text{in } \hat{\Omega}^e \tag{19}$$

where $\psi_i(\boldsymbol{\xi})$ are the $nd$-dimensional Lagrange interpolation functions, $c_i^e$ are the values of $c_{hp}(\mathbf{x})$ at the element nodes and $n = (p+1)^{nd}$ is the number of nodes in $\hat{\Omega}^e$. In addition, the quantity $p$ is the polynomial order of each interpolation function. There are a variety of ways in which high-order $nd$-dimensional interpolation functions may be formulated. For our analysis we construct these polynomial functions from tensor products of the one-dimensional $C^0$ spectral nodal interpolation functions

$$\varphi_j(\xi) = \frac{(\xi-1)(\xi+1) L'_p(\xi)}{p(p+1) L_p(\xi_j)(\xi - \xi_j)} \qquad \text{in } [-1, 1] \tag{20}$$

where $L_p(\xi)$ is the Legendre polynomial of order $p$. The quantities $\xi_j$ represent the locations of the nodes associated with the one-dimensional interpolants (with respect to the natural coordinate $\xi$). The one-dimensional nodal points are defined as the roots of the following expression

$$(\xi - 1)(\xi + 1) L'_p(\xi) = 0 \qquad \text{in } [-1, 1] \tag{21}$$

Figure 1 compares high-order interpolation functions for $p = 6$ using uniform nodes with corresponding interpolation functions using non-uniform nodes based on Gauss-Lobatto-Legendre points.

Multi-dimensional interpolation functions may be constructed from simple tensor products of the one-dimensional spectral interpolants. For example, in two-dimensions, the high-order interpolation functions may be defined as

$$\psi_i(\xi, \eta) = \varphi_j(\xi) \varphi_k(\eta) \qquad \text{in } \hat{\Omega}^e = [-1, 1]^2 \tag{22}$$

where $i = j + (k-1)(p+1)$ and $j, k = 1, \cdots, p+1$. Likewise, in three-dimensions, the interpolants are expressed as

$$\psi_i(\xi, \eta, \zeta) = \varphi_j(\xi) \varphi_k(\eta) \varphi_l(\zeta) \qquad \text{in } \hat{\Omega}^e = [-1, 1]^3 \tag{23}$$

where $i = j + [k - 1 + (l-1)(p+1)](p+1)$ and $j, k, l = 1, \cdots, p+1$. Finite elements constructed from tensor products of $\varphi_i(\xi)$ are commonly referred to as spectral elements in the literature [17]. Such elements are merely standard high order Lagrange type finite elements, where the locations of the unequally spaced nodes in $\hat{\Omega}^e$ are taken as tensor products of the roots of equation (21). Figure 2 illustrates various examples of high-order two-dimensional elements.

The finite element formulation naturally leads to a set of linear algebraic equations. The resulting coefficient matrix and force vector associated with the equations are constructed from the bilinear and linear forms appearing in equation (18) respectively. In this work we utilize the standard Gauss-Legendre quadrature rules for numerical integration of these quantities. We utilize *full integration* of all integrals and *do not* resort to selective under-integration of *any* terms in the coefficient matrix or force vector. All numerical results have been obtained using a quadrature rule of at least NGP $= p + 1$, where NGP represents the number of quadrature points in the direction of a given natural coordinate associated with $\hat{\Omega}^e$. For details on the computer implementation of the finite element method, including descriptions of the bijective isoparametric mapping $\bar{\Omega}^e \rightleftarrows \hat{\Omega}^e$ and the global assembly operator, we refer to the texts of Reddy [26] and Bathe [3]. For details on construction of the spectral interpolation functions, we refer to the book by Karniadakis and Sherwin [17].

## 5. REPRESENTATIVE NUMERICAL RESULTS

In this section, we illustrate the performance of high-order approximations under the classical single-field formulation and two least-squares-based formulations with respect to maximum principles and the non-negative constraint. Before we present numerical results, we briefly discuss how the numerical results under high-order approximations are visualized.

5.1. **Visualization of results using high-order approximations.** Once a numerical solution has been obtained using a particular finite element discretization, a given field variable may be evaluated at any point within a typical finite element using the standard interpolation formula given in equation (19). Use of this formula is crucial in evaluating the performance of high-order finite element models with respect to maximum principles and the non-negative constraint. Unfortunately, it is typically impossible to employ this scheme explicitly in the actual visualization of multi-dimensional numerical results obtained using high-order spectral/$hp$ finite element formulations. Visualization software (such as Tecplot [1]) allow for the post-processing of structured data associated with a given finite element mesh; however, such programs typically require data structures containing the element connectivity array of low-order elements only. To utilize standard visualization software, it therefore becomes necessary to convert data associated with a high-order spectral/$hp$ finite element mesh into data associated with a low-order finite element mesh.

Solution data associated with a high-order finite element mesh may be readily converted into solution data associated with a low-order finite element discretization through the creation of a set of fictitious low-order visualization elements, utilized for plotting purposes only. There are a variety of ways in which a low-order visualization mesh may be created. Perhaps, the simplest choice is to create a low-order mesh using only the actual nodes of the high-order finite element discretization. Unfortunately, the visualized numerical results will inevitably deviate from the actual finite element solution within a given element when such an approach is taken. In an effort to minimize visualization errors in the presentation of our numerical results, we utilize equation (19) to evaluate the numerical solution at a discrete number of points within a given element that is greater than the actual number of nodes of that element. For one-dimensional problems, we simply interpolate the numerical solution onto 100 grid points within each finite element. These grid points are then utilized explicitly to visualize the finite element solution. For two-dimensional problems, equation (19) is again employed to evaluate the numerical solution at 256 unequally spaced grid points within each element. These points are the effective nodes associated with a spectral/$hp$ finite element using a 15th order polynomial expansion (as geometrically characterized using the polynomial expansion associated with the actual element used in the numerical analysis). The refined set of points is then utilized to create a low-order visualization mesh that can be readily imported into Tecplot [1]. Alternatively, one can use any other procedure available in the literature for visualizing (scalar, vector and tensor) quantities interpolated using high-order approximations (for example, see [33, 27] and references therein).

We now illustrate the performance of $p$- and $h$-refinements on various one- and two-dimensional test problems. For two-dimensional problems, the regions of the violation of the non-negative constraint are indicated using white color.

5.2. **One-dimensional problem with zero forcing function.** This test problem is taken from Reference [22], which addressed the low-order approximation (i.e., $p = 1$). For completeness, we shall outline the problem. The domain is taken as $\Omega := (0, 1)$ with zero forcing function and non-zero Dirichlet boundary conditions. Mathematically, the test problem can be written as follows:

$$\alpha c - \frac{d^2 c}{dx^2} = 0 \quad \text{in } (0, 1) \tag{24a}$$

$$c(x = 0) = c(x = 1) = 1 \tag{24b}$$

where $\alpha$ is a non-negative constant. The analytical solution in terms of $\alpha$ is given by

$$c(x) = \frac{1 - \exp[-\sqrt{\alpha}]}{\exp[\sqrt{\alpha}] - \exp[-\sqrt{\alpha}]} \exp[\sqrt{\alpha}x] + \frac{\exp[\sqrt{\alpha}] - 1}{\exp[\sqrt{\alpha}] - \exp[-\sqrt{\alpha}]} \exp[-\sqrt{\alpha}x] \tag{25}$$

Herein we shall take $\alpha(x) = 1000$. It should be noted that the solution becomes steeper near the boundary as $\alpha$ becomes larger but the solution is still infinitely differentiable. The computational mesh is obtained by discretizing the domain using four equal-sized finite elements. The concentration profiles using the single-field formulation and the two least-squares-based formulations are, respectively, shown in Figures 3, 4 and 5. The flux profiles obtained using the LS2 formulation (which is a least-squares-based mixed formulation) is shown in Figure 6. In Reference [22], it has been shown that, for this test problem, the violation of the non-negative constraint vanishes with h-refinement under the classical single-field formulation. *From the numerical results presented in this subsection, one can conclude that the violation of the non-negative constraint also vanishes with p-refinement under the classical single-field and least-squares formulations for one-dimensional problems with zero forcing function.*

5.3. **One-dimensional problem with non-zero forcing function.** This test problem is taken from Reference [31]. Consider pure diffusion (that is, decay is neglected) in domain $\Omega := (-1, 1)$ with homogeneous Dirichlet boundary conditions. The forcing function is taken to be

$$f(\mathrm{x}) = 200 \exp[-10(\mathrm{x}+1)] \tag{26}$$

The test problem takes the following form:

$$-\frac{d^2 c}{d\mathrm{x}^2} = f(\mathrm{x}) \quad \text{in } \Omega \tag{27a}$$

$$c(\mathrm{x} = -1) = c(\mathrm{x} = 1) = 0 \tag{27b}$$

The analytical solution is given by

$$c(\mathrm{x}) = -2 \exp[-10(\mathrm{x}+1)] - (1 - \exp[-20])\mathrm{x} + (1 + \exp[-20]) \tag{28}$$

We shall mesh the computational domain using one element, and solve the problem for various polynomial approximations in the element. For this test problem, the LS1 and LS2 formulations are identical as $\alpha(\mathbf{x}) = 0$ and $\mathbf{D}(\mathbf{x}) = \mathbf{I}$. The obtained concentration profiles for various polynomial approximations under the classical single-field formulation and the least-squares formulations are shown in Figure 7. Although for this problem $f(\mathbf{x}) \geq 0 \; \forall \mathbf{x} \in \Omega$, the classical single-field formulation (for $p = 3$) and the least-squares formulation (for $p = 3$ and $p = 5$) violate the non-negative constraint. The reason for the violation is that, under the classical single-field formulation, the $L_2$ projection of the forcing function onto the space spanned by $p = 3$ has negative values near $\mathrm{x} = 1$. A similar reason holds for the least-squares formulation.

5.4. **Two-dimensional isotropic diffusion.** This test problem was proposed by Burman and Ern [6]. The computational domain is a rectangle $\Omega := (0, 1) \times (0, 0.3)$ with homogeneous Dirichlet boundary conditions. The decay coefficient $\alpha(\mathbf{x})$ is taken to be zero. The diffusivity tensor is the second-order identity tensor. That is,

$$\mathbf{D}(\mathbf{x}) = \mathbf{I} \tag{29}$$

The forcing function is taken as follows:

$$f(\mathbf{x}) = \begin{cases} 1 & \text{if } (\mathrm{x}, \mathrm{y}) \in [0, 0.5] \times [0, 0.075] \\ 0 & \text{otherwise} \end{cases} \tag{30}$$

TABLE 1. Two-dimensional isotropic diffusion: This table shows the variation of the minimum concentration with respect to $p$-refinement and $h$-refinement. Since the medium is isotropic, the violation of the non-negative constraint vanishes with $h$-refinement. Although the forcing function is non-negative, the violation of the non-negative constraint did not vanish with $p$-refinement. One of the reasons is that the $L_2$ projection of the forcing function onto the polynomial space is *not* non-negative.

| $h/p$ order | $p$-refinement | | $h$-refinement | |
| --- | --- | --- | --- | --- |
| | Single-field | Least-squares | Single-field | Least-squares |
| 1 | 0 | 0 | 0 | 0 |
| 2 | $-1.35 \times 10^{-4}$ | $-1.41 \times 10^{-4}$ | $-8.84 \times 10^{-5}$ | $-1.08 \times 10^{-4}$ |
| 3 | $-3.13 \times 10^{-5}$ | $-3.17 \times 10^{-5}$ | $-1.81 \times 10^{-5}$ | $-2.57 \times 10^{-5}$ |
| 4 | $-2.04 \times 10^{-5}$ | $-2.08 \times 10^{-5}$ | 0 | 0 |
| 5 | $-4.42 \times 10^{-6}$ | $-4.46 \times 10^{-6}$ | 0 | 0 |
| 6 | $-2.77 \times 10^{-6}$ | $-2.79 \times 10^{-6}$ | 0 | 0 |
| 7 | $-1.99 \times 10^{-7}$ | $-2.01 \times 10^{-7}$ | 0 | 0 |
| 8 | $-6.22 \times 10^{-7}$ | $-6.24 \times 10^{-7}$ | 0 | 0 |
| 9 | $-2.60 \times 10^{-8}$ | $-2.60 \times 10^{-8}$ | 0 | 0 |
| 10 | $-2.06 \times 10^{-7}$ | $-2.06 \times 10^{-7}$ | 0 | 0 |

Since the forcing function $f(\mathbf{x}) \geq 0$ in $\Omega$ and $c^{\mathrm{p}}(\mathbf{x}) = 0$ on $\partial \Omega$, from the maximum principle given in Theorem 1 we have $c(\mathbf{x}) \geq 0$ in $\bar{\Omega}$. It should be noted that since the diffusivity tensor is isotropic and there is no decay, both the least-squares formulations are identical (see equation (13)). The computational and visualization meshes are shown in Figure 8. The concentration profiles obtained under the single-field and the least-squares formulations for various polynomial approximations are, respectively, shown in Figures 9 and 10.

The minimum concentration under the single-field and the least-squares formulations for various order of $p$- and $h$-refinements are given in Table 1. As one can see, there is violation of the non-negative constraint under $p$-refinement, and the violation vanishes under $h$-refinement. As discussed towards the end of Section 2, the reason is that the $L_2$ projection of the forcing function on the polynomial space need not be non-negative under high-order approximations (even though the forcing function is non-negative). Note that the $L_2$ projection of the forcing function onto $p = 1$ polynomial space is non-negative.

5.5. **Non-uniform anisotropic media.** This test problem was originally proposed by Le Potier [24], and has been employed in many other numerical studies on maximum principles and the non-negative constraint using the low-order approximation (e.g., see References [21, 23]). The test problem assumes $\alpha(\mathbf{x}) = 0$. The computational domain is a square $\Omega := (0, 0.5) \times (0, 0.5)$ with homogeneous Dirichlet boundary conditions enforced on the whole boundary. The diffusivity tensor is given by

$$\mathbf{D}(\mathbf{x}) = \begin{pmatrix} y^2 + \epsilon x^2 & -(1-\epsilon)xy \\ -(1-\epsilon)xy & x^2 + \epsilon y^2 \end{pmatrix} \tag{31}$$

with $\epsilon = 10^{-3}$. The forcing function is given by

$$f(\mathbf{x}) = \begin{cases} 1 & (x, y) \in [0.125, 0.375] \times [0.125, 0.375] \\ 0 & \text{otherwise} \end{cases} \tag{32}$$

Two different meshes (one is a structured mesh and the other is an unstructured mesh) are employed, and are shown in Figures 11 and 12. The performance of the single-field, LS1 and LS2 formulations using the structured mesh under both $p$- and $h$-refinements is illustrated in Figures 13, 14 and 15. The performance of these formulations using the the unstructured mesh is illustrated in Figures 16, 17 and 18. The minimum concentration under these three formulations is shown in Figure 19. Since the anisotropy is strong, there will be violation of the non-negative constraint even on structured mesh for both $p$- and $h$-refinements. For this test problem, it should be noted that LS1 formulation did not perform as well as LS2 formulation. This is expected for problems involving heterogeneous and anisotropic media as the weight(s) in defining least-squares functional play a crucial role in the accuracy of the numerical results. (See equation (13) to note the different weights used in constructing LS1 and LS2 formulations.)

5.6. **Anisotropic diffusion in a square domain with a hole.** This problem has been used in References [21, 19, 23, 22] with respect to the enforcement of the non-negative constraint but in the context of low-order approximation. The computational domain is a bi-unit square with a square hole of dimension $[4/9, 5/9] \times [4/9, 5/9]$. The forcing function is taken to be $f(\mathbf{x}) = 0$, and the decay coefficient is $\alpha(\mathbf{x}) = 0$. On the external boundary $c^{\mathrm{p}}(\mathbf{x}) = 0$ is prescribed, and on the internal boundary $c^{\mathrm{p}}(\mathbf{x}) = 2$ is prescribed. The diffusivity tensor is given by

$$\mathbf{D}(\mathbf{x}) = \begin{pmatrix} \cos(\theta) & \sin(\theta) \\ -\sin(\theta) & \cos(\theta) \end{pmatrix} \begin{pmatrix} k_1 & 0 \\ 0 & k_2 \end{pmatrix} \begin{pmatrix} \cos(\theta) & -\sin(\theta) \\ \sin(\theta) & \cos(\theta) \end{pmatrix} \tag{33}$$

Herein, we have taken $\theta = \pi/6$, and considered two different sets for diffusivity coefficients: $(k_1, k_2) = (1, 100)$ and $(k_1, k_2) = (1, 10000)$. The computational and visualization meshes are shown in Figure 20. The concentration profiles under the single-field, LS1 and LS2 formulations for various polynomial approximations are, respectively, shown in Figures 21, 22 and 23. Figure 24 shows the minimum concentration under these three formulations for both sets of diffusivity coefficients. From this figure, it is evident that the greater is disparity between the diffusivity coefficients $k_1$ and $k_2$ the greater is the violation of the non-negative constraint. Moreover, the extent of the violation did not decrease with with $p$- and $h$-refinements for strong anisotropic medium (see the case $k_1 = 1$ and $k_2 = 10000$). Hence, one can conclude that for strong anisotropic medium, $p$- and $h$-refinement do not eliminate the violation of the non-negative constraint.

## 6. CONCLUSIONS

In this paper, we considered the diffusion-type equation, which is a second-order elliptic partial differential equation. This particular equation is known to satisfy maximum principles and the non-negative constraint under certain conditions. Many popular numerical formulations do not satisfy either maximum principles or the non-negative constraint. In the literature, it has been documented the performance of finite element formulations with respect to maximum principles for low-order elements. Herein, we considered the classical single-field formulation as well as two least-squares-based mixed formulations. We have systematically documented the performance of these formulations with respect to maximum principles and the non-negative constraint under high-order approximations. The main findings of this paper can be summarized as follows:

(a) For one-dimensional problems, it is well-known that a uniform mesh is sufficient to satisfy the non-negative constraint and maximum principles under the low-order approximation. (For non-zero decay coefficient, the size of the element should be smaller than a critical size.) However, a uniform mesh is not sufficient for high-order approximations. (See subsections 5.2 and 5.3.)

(b) For isotropic two-dimensional problems, a well-centered triangular mesh or a mesh with rectangular elements with aspect ratio between $1/\sqrt{2}$ and $\sqrt{2}$ is sufficient for the low-order approximation to satisfy the non-negative constraint and maximum principles. Again, these conditions are not sufficient for high-order approximations. (See subsection 5.4.)

(c) For anisotropic two-dimensional problems, any restrictions on the mesh will not be sufficient, in general, to satisfy the non-negative constraint and maximum principles for both low-order and high-order approximations. (See subsections 5.5 and 5.6.)

(d) The performance of a least-squares formulation depends on the weight(s) used in defining the least-squares functional. (See subsection 5.5, and Figures 14 and 18.)

We shall conclude this paper with the following statement with a hope that it will motivate applied mathematicians, computational mechanicians, and numerical analysts to work on an interesting problem: *A finite-element-based formulation or methodology that satisfies the non-negative constraint and maximum principles for anisotropic diffusion on general computational grids under* *high-order approximations is currently an unsolved problem with many important applications in engineering and applied sciences.*

## ACKNOWLEDGMENTS

The first author (G. S. Payette) acknowledges the Sandia National Laboratories/Texas A&M University Excellence in Engineering Graduate Research Fellowship. The second author (K. B. Nakshatrala) acknowledges the financial support from Sandia National Laboratories through the Laboratory Directed Research and Development program. Sandia is a multiprogram laboratory operated by Sandia Corporation, a Lockheed Martin Company, for the United States Department of Energy's National Nuclear Security Administration under Contract DE-AC04-94AL85000. The third author (J. N. Reddy) gratefully acknowledges the support by the Oscar S. Wyatt Endowed Chair. The opinions expressed in this paper are those of the authors and do not necessarily reflect that of the sponsors.

## References

[1] *Tecplot 360: User's Manual*. URL: http://www.tecplot.com, Bellevue, Washington, USA, 2008.

[2] R. B. Bapat and T. E. S. Raghavan. *Non-negative Matrices and Applications*. Cambridge University Press, Cambridge, UK, 1997.

[3] K.-J. Bathe. *Finite element procedures*. Prentice Hall, Upper Saddle River, New Jersey, 1996.

[4] P. B. Bochev and M. D. Gunzburger. *Least-Squares Finite Element Methods*. Springer, New York, USA, 2009.

[5] M. Borsuk and V. Kondratiev. *Elliptic Boundary Value Problems of Second Order in Piecewise Smooth Domains*. Elsevier Science, San Diego, USA, 2006.

[6] E. Burman and A. Ern. Discrete maximum principle for Galerkin approximations of the Laplace operator on arbitrary meshes. *Comptes Rendum Mathematique*, 338:641–646, 2004.

[7] Y.-Z. Chen and L.-C. Wu. *Second Order Elliptic Equations and Elliptic Systems*. American Mathematical Society, New York, USA, 2004.

[8] I. Christie and C. Hall. The maximum principle for bilinear elements. *International Journal for Numerical Methods in Engineering*, 20:549–553, 1984.

[9] P. G. Ciarlet and P-A. Raviart. Maximum principle and uniform convergence for the finite element method. *Computer Methods in Applied Methods and Engineering*, 2:17–31, 1973.

[10] D. Gilbarg and N. S. Trudinger. *Elliptic Partial Differential Equations of Second Order*. Springer, New York, USA, 2001.

[11] M. E. Gurtin. *An Introduction to Continuum Mechanics*. Academic Press, San Diego, USA, 1981.

[12] Q. Han and F. Lin. *Elliptic Partial Differential Equations*. American Mathematical Society, Providence, Rhode Island, USA, 2000.

[13] P. Herrera and A. Valocchi. Positive solution of two-dimensional solute transport in heterogeneous aquifers. *Ground Water*, 44:803–813, 2006.

[14] D. W. Höhn and D. H. D. Mittelmann. Some remarks on the discrete maximum-principle for finite elements in higher order. *Computing*, 27:145–154, 1981.

[15] E. Hopf. Elementäre Bemerkungen über die Lösungen partieller Differentialgleichungen zweiter Ordnung vom elliptischen Typus. *Sitzungsberichte der Berliner Akademie der Wissenschaften*, 19:147–152, 1927.

[16] B. Jiang. *The Least-Squares Finite Element Method: Theory and Applications in Computational Fluid Dynamics and Electromagnetics*. Springer-Verlag, New York, USA, 1998.

[17] G. E. Karniadakis and S. J. Sherwin. *Spectral/hp Element Methods for CFD*. Oxford University Press, Oxford, 1999.

[18] D. Kuzmin and S. Turek. Flux correction tools for finite elements. *Journal of Computational Physics*, 175:525–558, 2002.

[19] K. Lipnikov, M. Shashkov, D. Svyatskiy, and Y. Vassilevski. Monotone finite volume schemes for diffusion equations on unstructured triangular and shape-regular polygonal meshes. *Journal of Computational Physics*, 227:492–512, 2007.

[20] K. Lipnikov, D. Svyatskiy, and Y. Vassilevski. Interpolation-free monotone finite volume method for diffusion equations on polygonal meshes. *Journal of Computational Physics*, 228:703–716, 2009.

[21] R. Liska and M. Shashkov. Enforcing the discrete maximum principle for linear finite element solutions for elliptic problems. *Communications in Computational Physics*, 3:852–877, 2008.

[22] H. Nagarajan and K. B. Nakshatrala. Enforcing the non-negativity constraint and maximum principles for diffusion with decay on general computational grids. *International Journal for Numerical Methods in Fluids*, In print, DOI: 10.1002/fld.2389, 2010.

[23] K. B. Nakshatrala and A. J. Valocchi. Non-negative mixed finite element formulations for a tensorial diffusion equation. *Journal of Computational Physics*, 228:6726–6752, 2009.

[24] C. Le Potier. Finite volume monotone scheme for highly anisotropic diffusion operators on unstructured triangular meshes. *Comptes Rendus Mathematique*, 341:787–792, 2005.

[25] P. Pucci and J. Serrin. *The Maximum Principle*. Birkhäuser Verlag, Basel, Switzerland, 2007.

[26] J. N. Reddy. *An Introduction to the Finite Element Method*. McGraw-Hill, New York, 3rd edition, 2006.

[27] J.-F. Remavle, N. Chevaugeon, E. Marchandise, and C. Geuzaine. Efficient visualization of high-order finite elements. *International Journal for Numerical Methods in Engineering*, 69:750–771, 2007.

[28] Y. Saad. *Iterative Methods for Sparse Linear Systems*. SIAM, Philadelphia, USA, 2003.

[29] J. Serrin. Commentary on the Hopf strong maximum principle. In C. S. Morawetz, J. Serrin, and Y. G. Sinai, editors, *Selected Works of Eberhard Hopf with Commentaries*, pages 9–14. American Mathematical Society, Providence, 2002.

[30] P. Solin and T. Vejchodsky. A weak discrete maximum principle for *hp*-FEM. *Journal of Computational and Applied Mathematics*, 209:54–65, 2007.

[31] P. Solin, T. Vejchodsky, and R. Araiza. Discrete conservation of nonnegativity for elliptic problems solved by *hp*-FEM. *Mathematics and Computers in Simulations*, 76:205–210, 2007.

[32] M. E. Taylor. *Partial Differential Equations I: Basic Theory*. Springer-Verlag, New York, USA, 1996.

[33] D. C. Thompson and P. P. Pébay. Visualizing higher order finite elements: Final report. Technical Report SAND2005-6999, Sandia National Laboratories, November 2005.

[34] R. Varga. On discrete maximum principle. *SIAM Journal on Numerical Analysis*, 3:355–359, 1966.

[35] T. Vejchodsky and P. Solin. Discrete maximum principle for high-order finite elements in 1D. *Mathematics of Computation*, 76:1833–1846, 2007.

[36] E. G. Yanik. Sufficient conditions for a discrete maximum principle for high order collocation methods. *Computers and Mathematics with Applications*, 17:1431–1434, 1989.

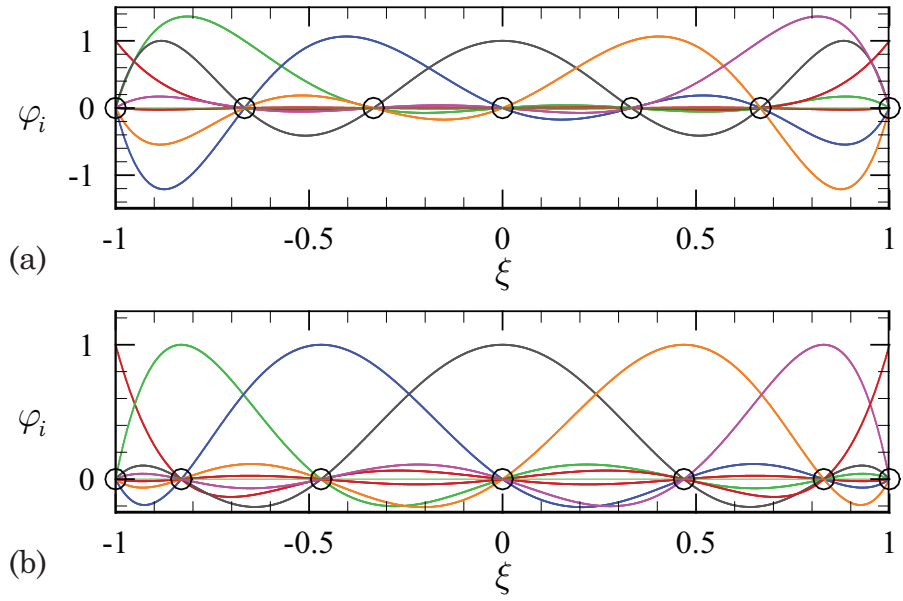

FIGURE 1. High polynomial order one-dimensional $C^0$ Lagrange interpolation functions (cases shown are for $p = 6$) with (a) equal spacing of the element nodes, and (b) unequal nodal spacing associated with Gauss-Lobatto-Legendre points.

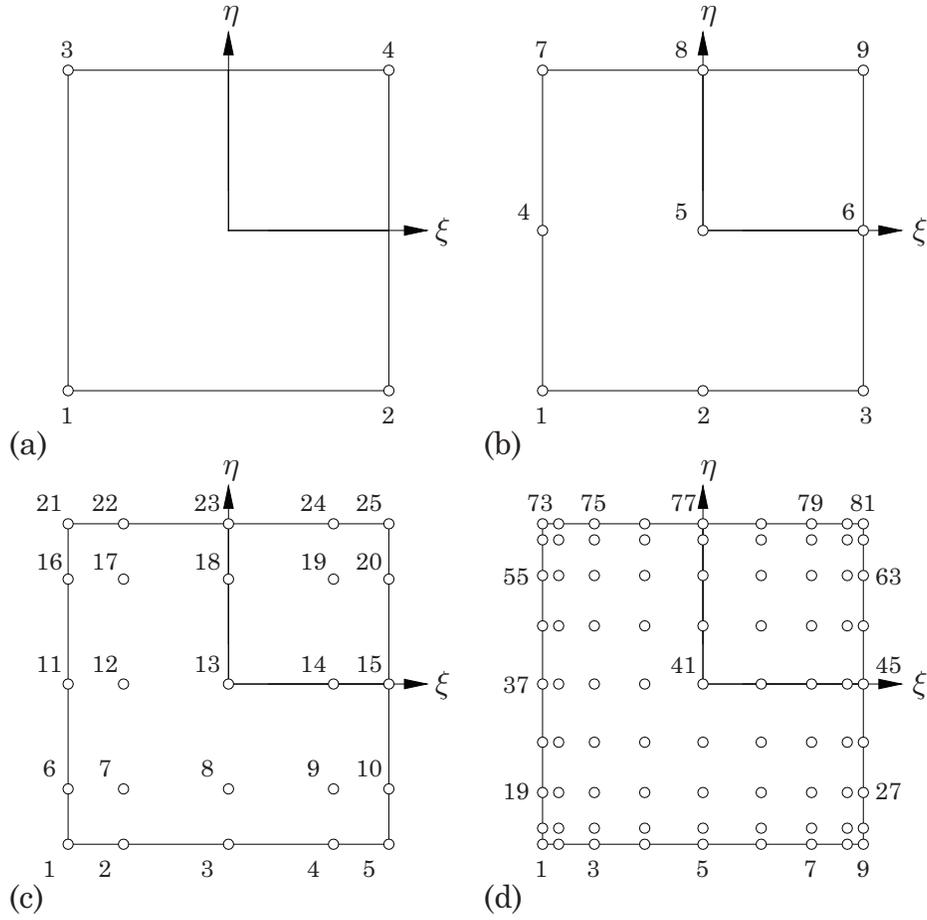

FIGURE 2. Examples of various high-order two-dimensional master elements $\hat{\Omega}^e$ used in the present formulation: (a) four-node element ($p = 1$), (b) nine-node element ($p = 2$), (c) twenty-five-node element ($p = 4$), and (d) eighty-one-node element ($p = 8$).

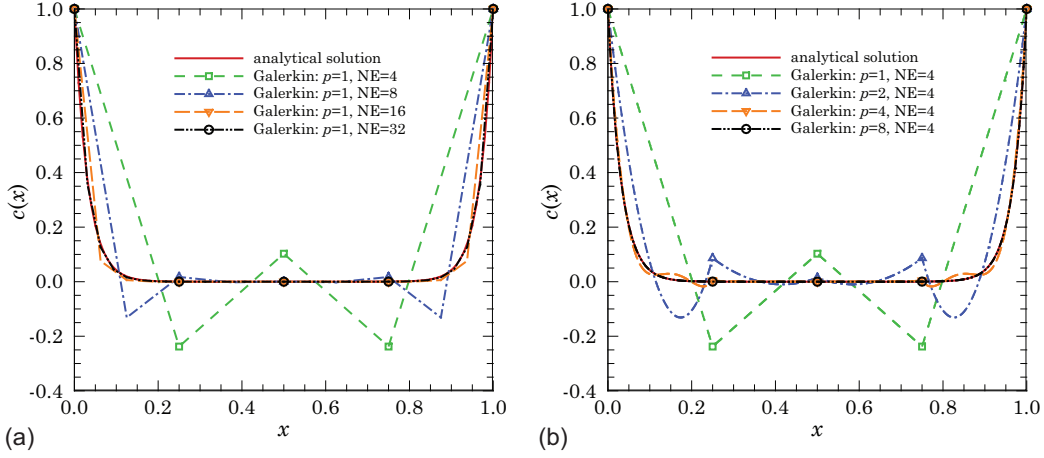

FIGURE 3. One-dimensional problem with zero forcing function: The figure illustrates the performance of the classical single-field formulation (which is also referred to as the Galerkin formulation) with respect to the non-negative constraint. In the right figure, we show the classical single-field formulation with respect to $p$-refinement. For comparison, the left figure shows the performance for an equivalent $h$-refinement. For this test problem, the violation of the non-negative constraint decreases with both $h$- and $p$-refinements under the classical single-field formulation.

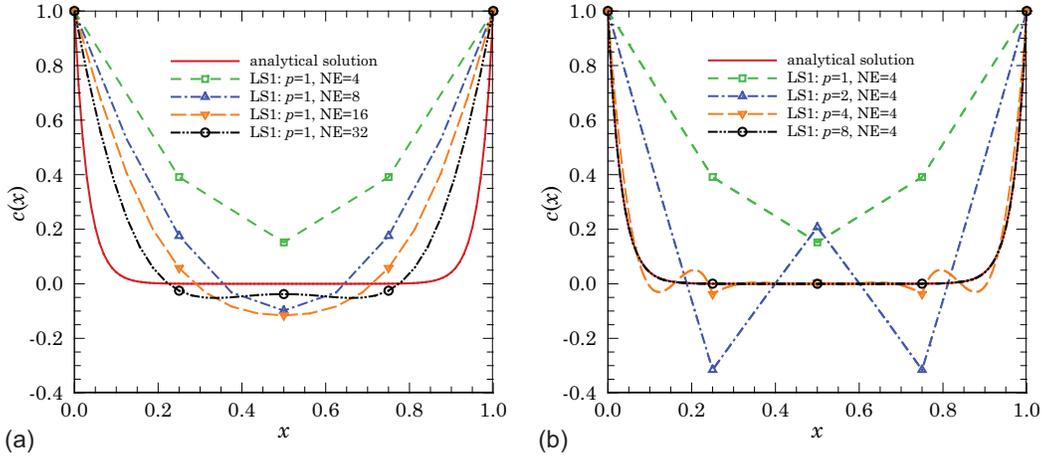

FIGURE 4. One-dimensional problem with zero forcing function: The figure illustrates the performance of the first-form of the least-squares formulation (LS1 formulation) with respect to the non-negative constraint. In the right figure, we show the performance of the formulation with respect to $p$-refinement. For comparison, the left figure shows the performance for an equivalent $h$-refinement. For this problem, the violation of the non-negative constraint decreases with both $h$- and $p$-refinements under the LS1 formulation.

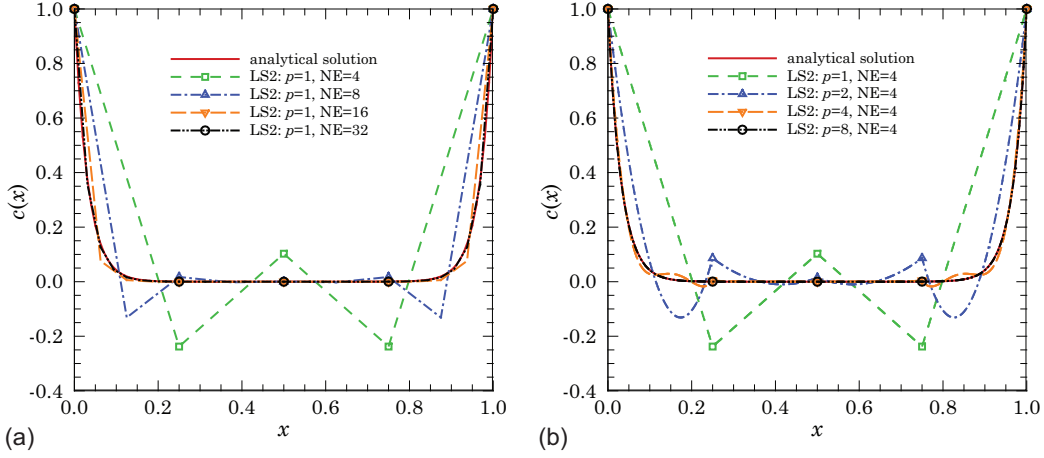

FIGURE 5. One-dimensional problem with zero forcing function: The figure illustrates the performance of the second-form of the least-squares formulation (LS2 formulation). In the right figure, we show the performance of the formulation due to $p$-refinement. For comparison, the left figure shows the performance for equivalent $h$-refinement. For this problem, the violation of the non-negative constraint decreases with both $h$- and $p$-refinements under the LS1 formulation.

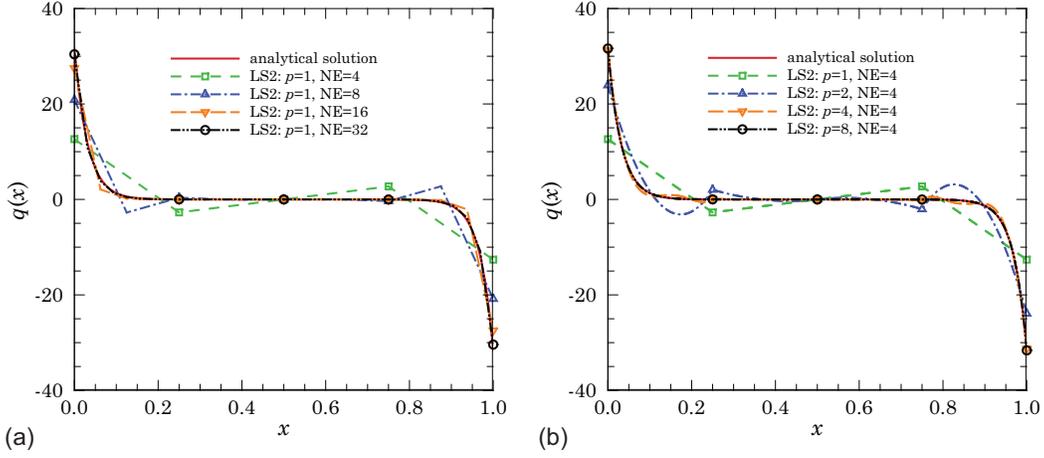

FIGURE 6. One-dimensional problem with zero forcing function: The figure shows the flux obtained using the second-form of the least-squares formulation (LS2 formulation). In the right figure, we show the performance of the formulation due to $p$-refinement. For comparison, the left figure shows the performance for equivalent $h$-refinement.

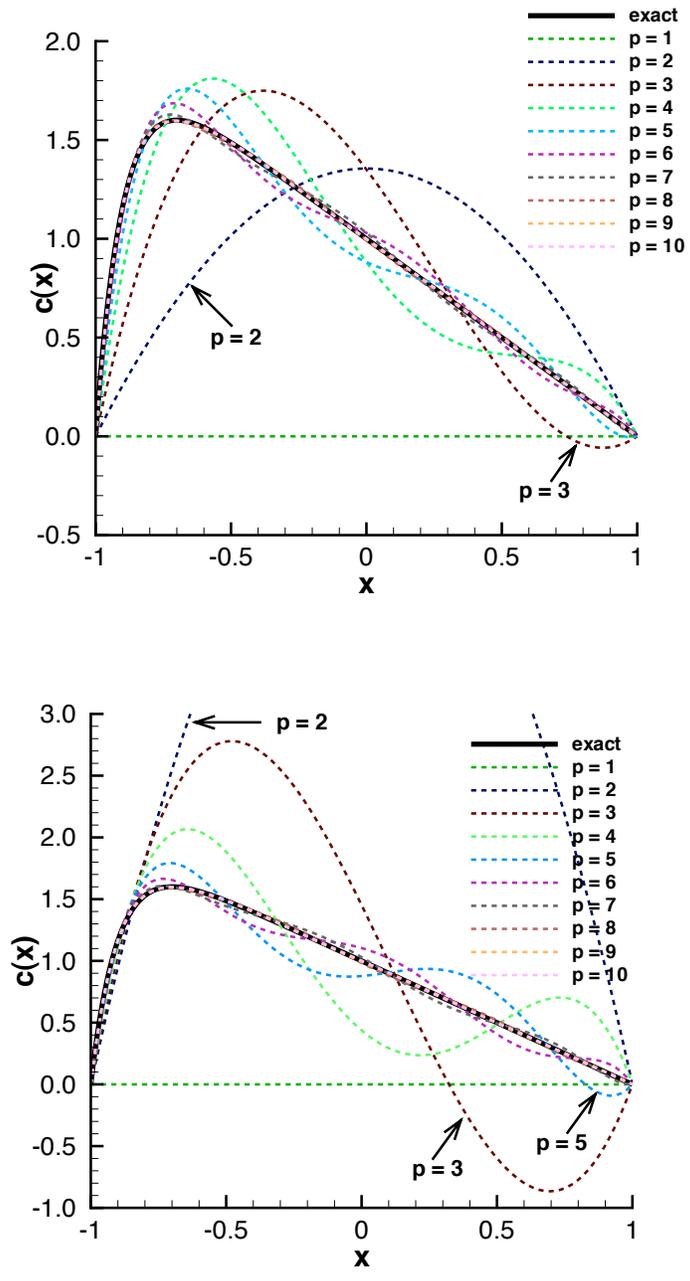

FIGURE 7. One-dimensional problem with non-zero forcing function: This figure shows the obtained concentration profiles for various orders of polynomial approximations under the classical single-field formulation (top figure) and the least-squares formulation (bottom figure). Although $f(\mathbf{x}) \geq 0 \ \forall \mathbf{x} \in \Omega$, both the classical single-field formulation (for $p = 3$) and the least-squares formulation (for $p = 3$ and $p = 5$) violate the non-negative constraint.

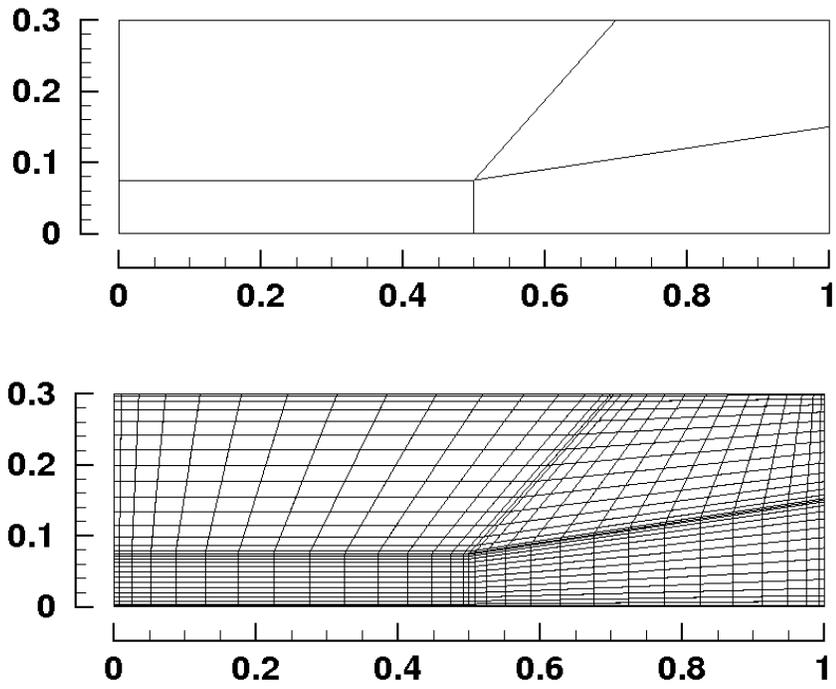

FIGURE 8. Two-dimensional isotropic diffusion: The top figure shows the mesh used in the numerical simulation, and the bottom figure shows the corresponding mesh used for visualization under high-order approximations. For more details on visualizing results using high-order approximations see subsection 5.1.

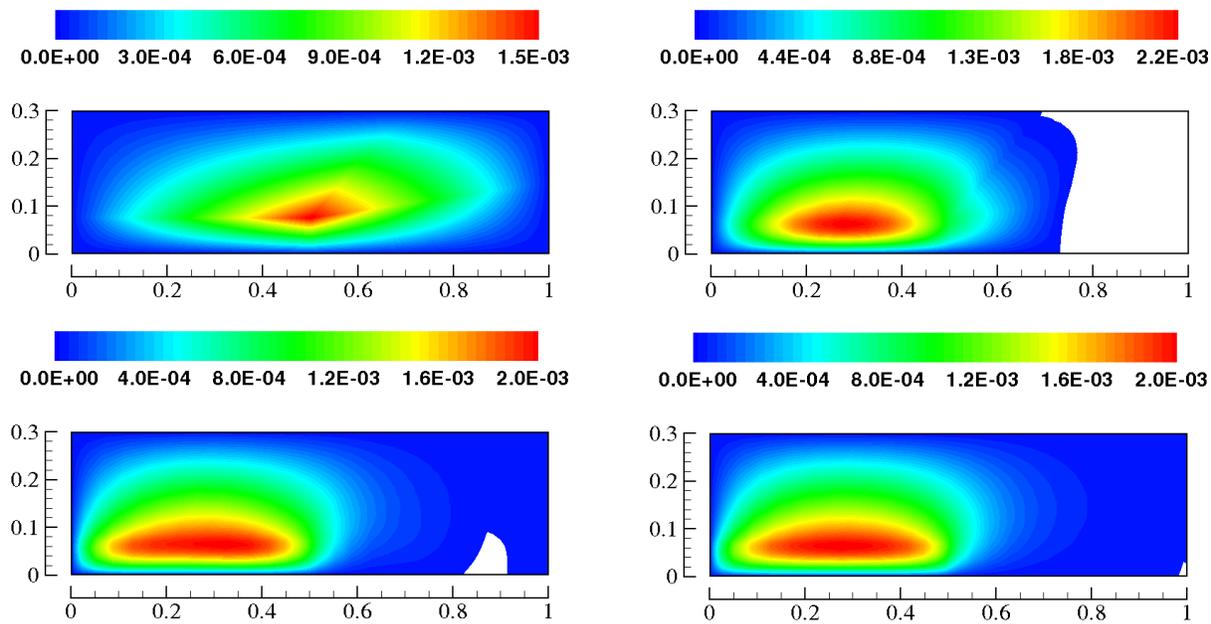

Figure 9. Two-dimensional isotropic diffusion: This figure shows the contours of the concentration obtained using the *single-field formulation* for various polynomial approximations. The top figures are for $p = 1$ and $p = 2$, and the bottom figures are for $p = 5$ and $p = 10$. The regions in which the non-negative constraint is violated is shown in white color.

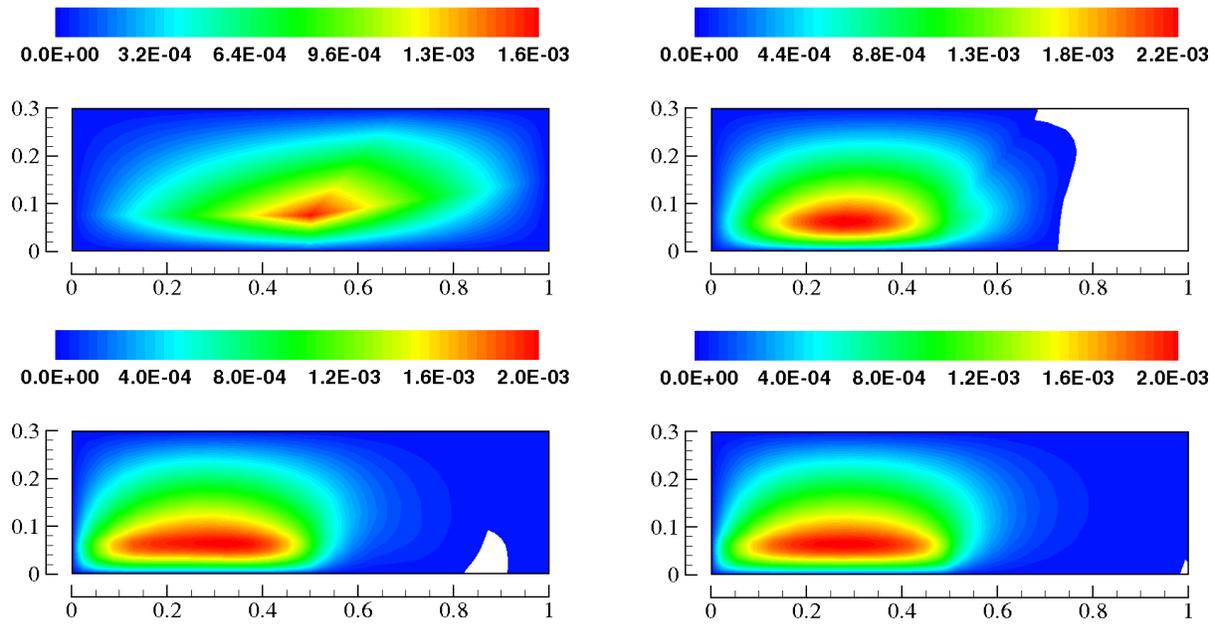

FIGURE 10. Two-dimensional isotropic diffusion: This figure shows the contours of the concentration obtained using the *least-squares formulation* for various polynomial approximations. The top figures are for $p = 1$ and $p = 2$, and the bottom figures are for $p = 5$ and $p = 10$. (Since the diffusivity is isotropic and decay coefficient is zero, both LS1 and LS2 formulations are identical.) The regions in which the non-negative constraint is violated is shown in white color.

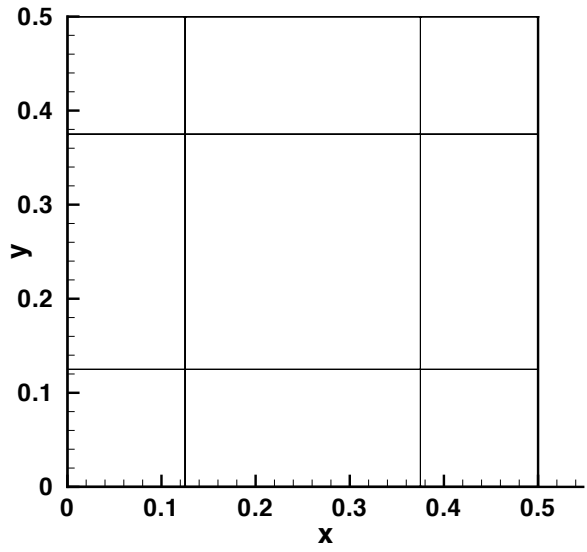

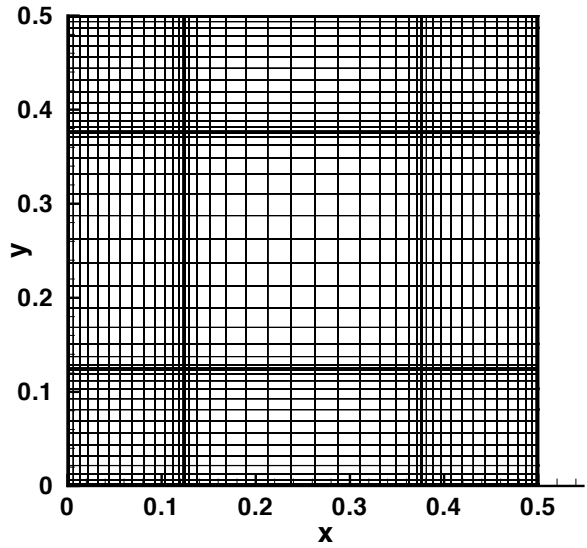

FIGURE 11. Non-uniform anisotropic media: The top figure shows the *structured* mesh used in the numerical simulation, and the bottom figure shows the corresponding mesh used for visualization under high-order approximations. For more details on visualizing results using high-order approximations see subsection 5.1.

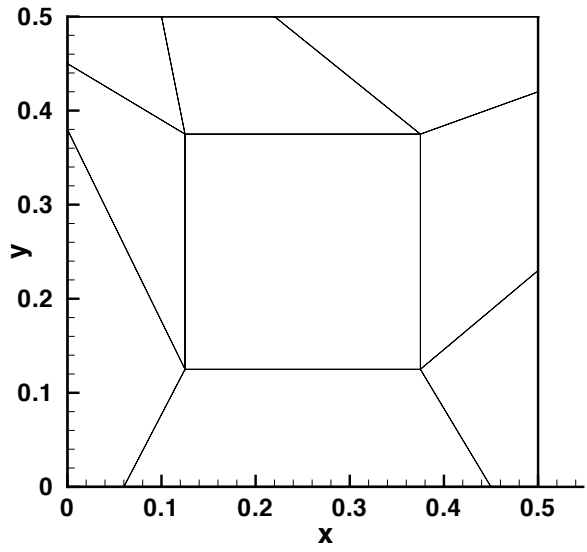

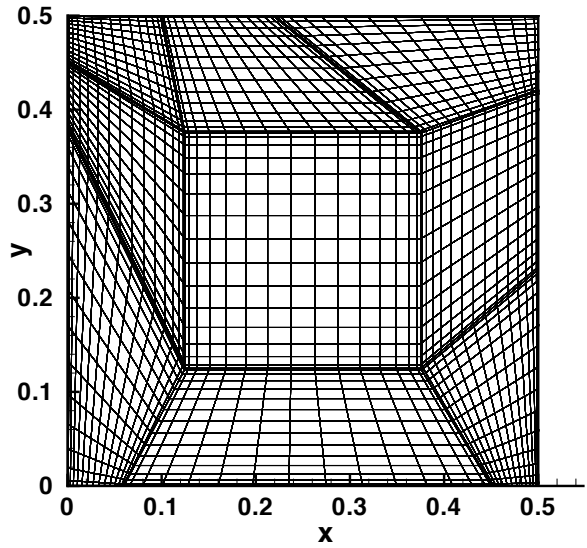

FIGURE 12. Non-uniform anisotropic media: The top figure shows the *unstructured* mesh used in the numerical simulation, and the bottom figure shows the corresponding mesh used for visualization under high-order approximations. For more details on visualizing results using high-order approximations see subsection 5.1.

Gregory S. Payette, Graduate student, Department of Mechanical Engineering, Texas A&M University, College Station, Texas - 77843.

*E-mail address*: `paye4446@tamu.edu`

Correspondence to: Dr. Kalyana Babu Nakshatrala, Department of Mechanical Engineering, 216 Engineering/Physics Building, Texas A&M University, College Station, Texas - 77843. TEL: +1-979-845-1292

*E-mail address*: `knakshatrala@tamu.edu`

Dr. J. N. Reddy, Department of Mechanical Engineering, Texas A&M University, 210 Engineering/Physics Building, Texas A&M University, College Station, Texas - 77843. TEL:+1-979-862-2417

*E-mail address*: `jnredddy@tamu.edu`

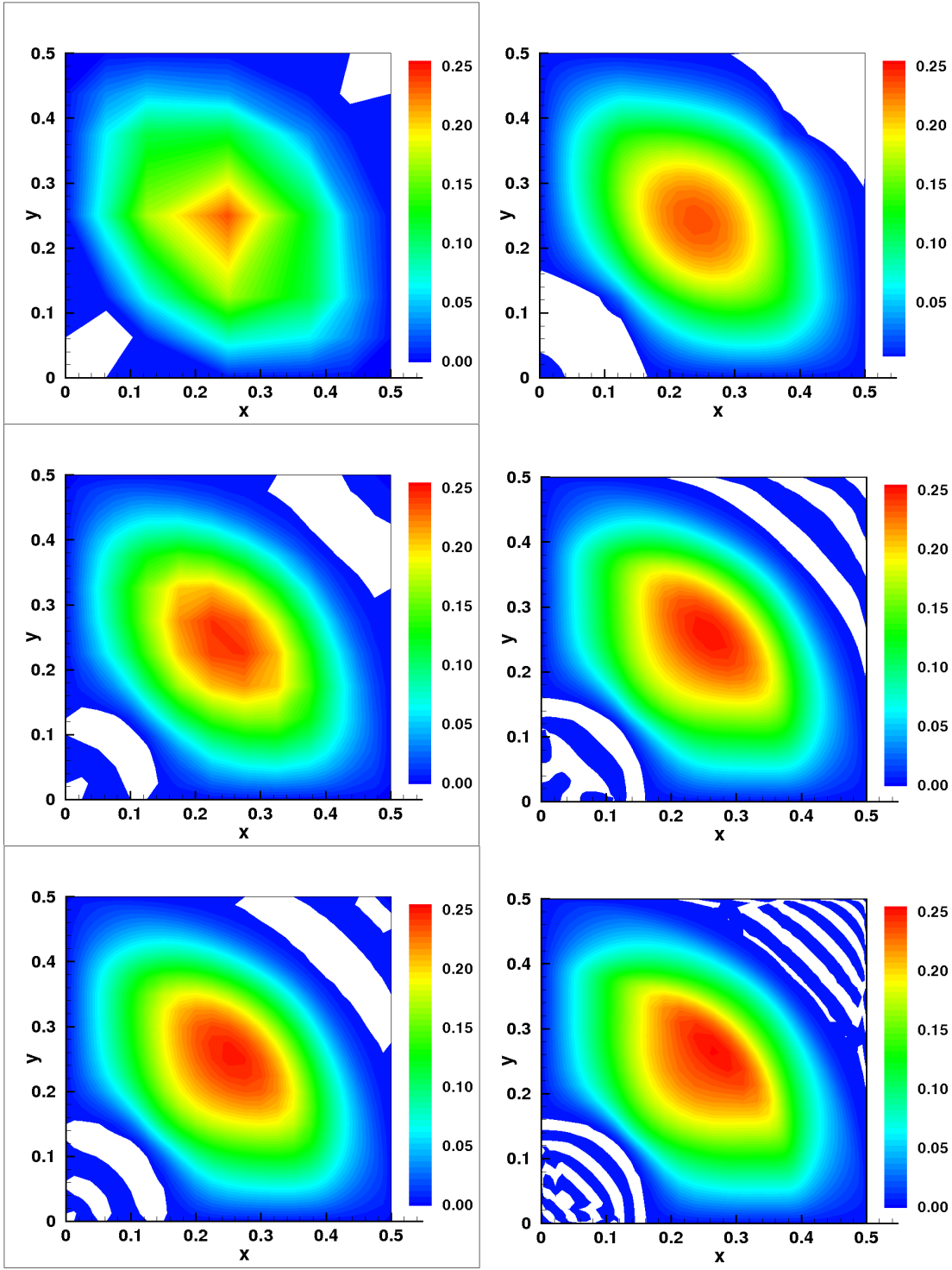

FIGURE 13. Non-uniform anisotropic media: This figure shows the concentration contours obtained using the *single-field formulation* and structured mesh (which is show in Figure 11). The left figures are using $h$-refinement, and the right figures are using $p$-refinement. The top figures are for $h/p = 2$, the middle figures are for $h/p = 5$, and the bottom two figures are for $h/p = 10$.

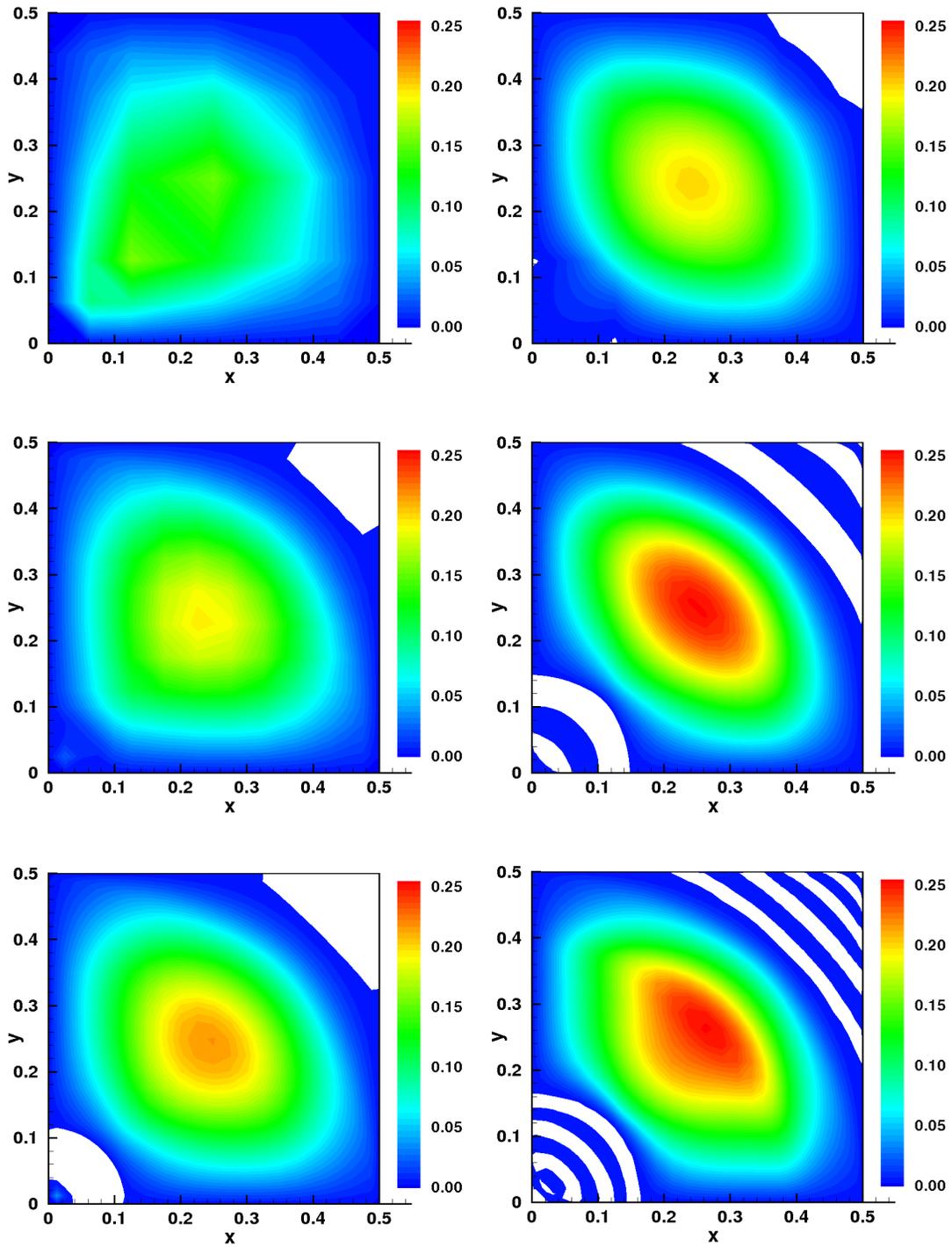

FIGURE 14. Non-uniform anisotropic media: This figure shows the concentration contours obtained using the *LS1 formulation* and *structured mesh* (which is shown in Figure 11). The left figures are using $h$-refinement, and the right figures are using $p$-refinement. The top figures are for $h/p = 2$, the middle figures are for $h/p = 5$, and the bottom two figures are for $h/p = 10$.

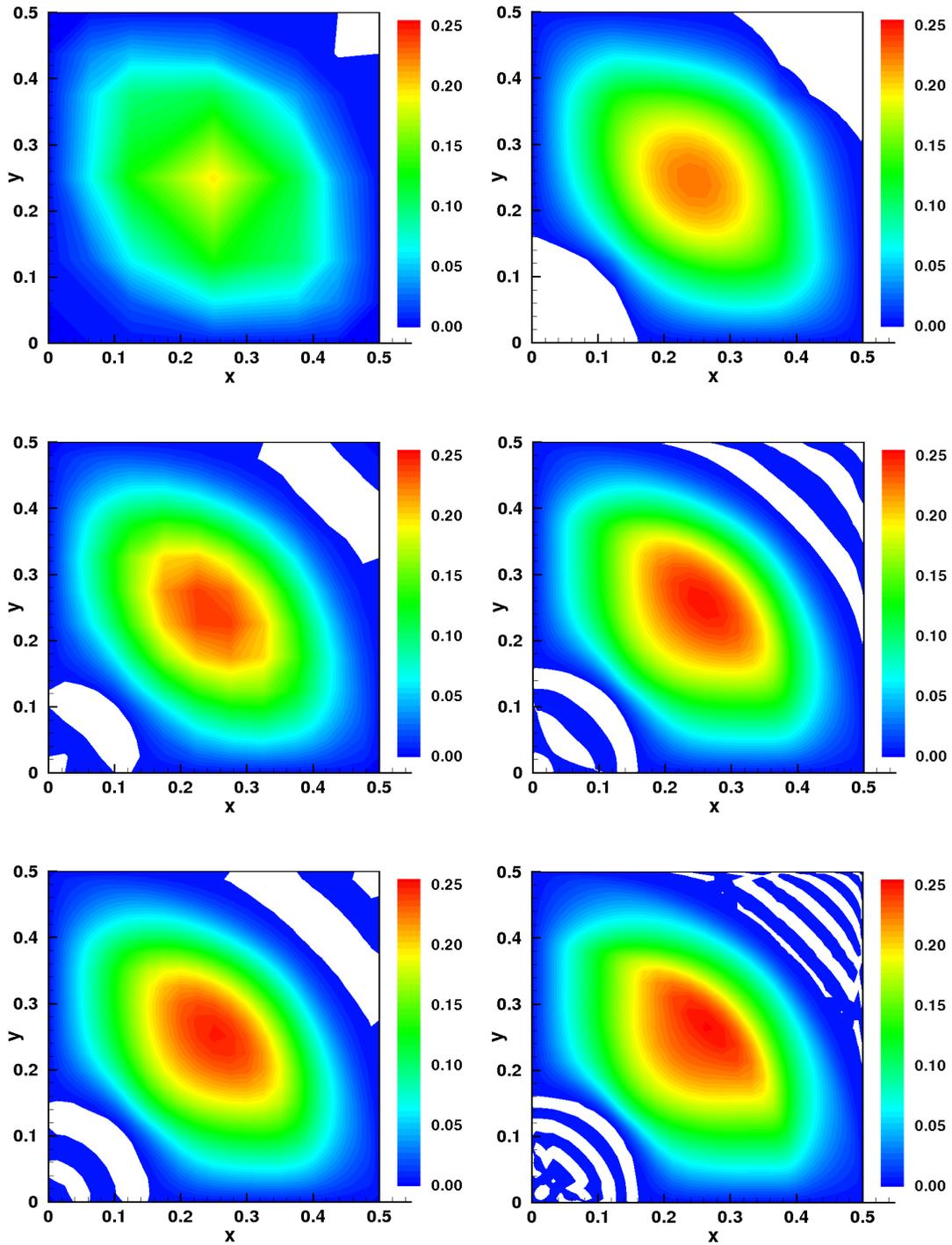

FIGURE 15. Non-uniform anisotropic media: This figure shows the concentration contours obtained using the *LS2 formulation* and *structured mesh* (which is shown in Figure 11). The left figures are using $h$-refinement, and the right figures are using $p$-refinement. The top figures are for $h/p = 2$, the middle figures are for $h/p = 5$, and the bottom two figures are for $h/p =$ 10.

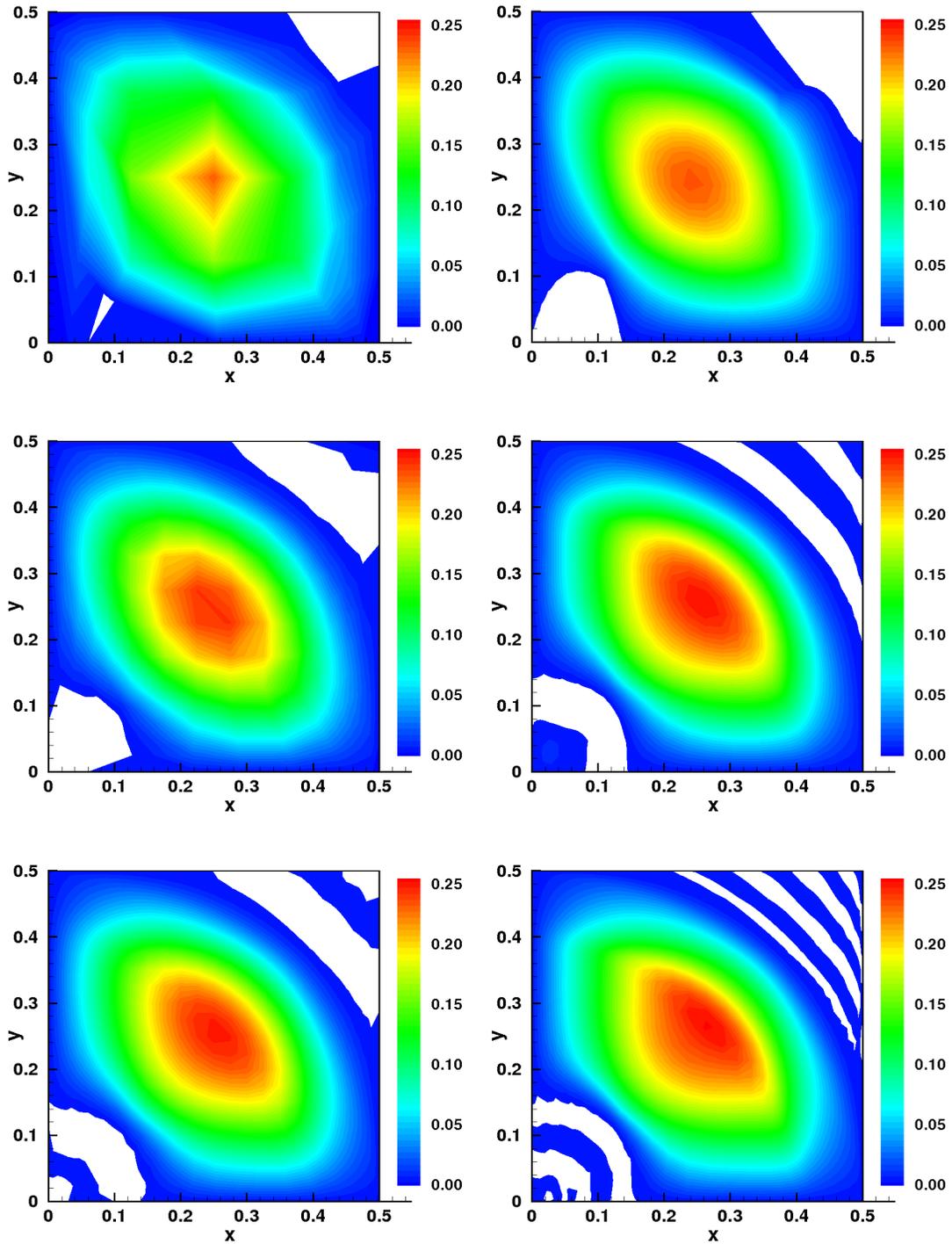

FIGURE 16. Non-uniform anisotropic media: This figure shows the concentration contours obtained using the *single-field formulation* and *unstructured mesh* (which is shown in Figure 12). The left figures are using $h$-refinement, and the right figures are using $p$-refinement. The top figures are for $h/p = 2$, the middle figures are for $h/p = 5$, and the bottom two figures are for $h/p = 10$.

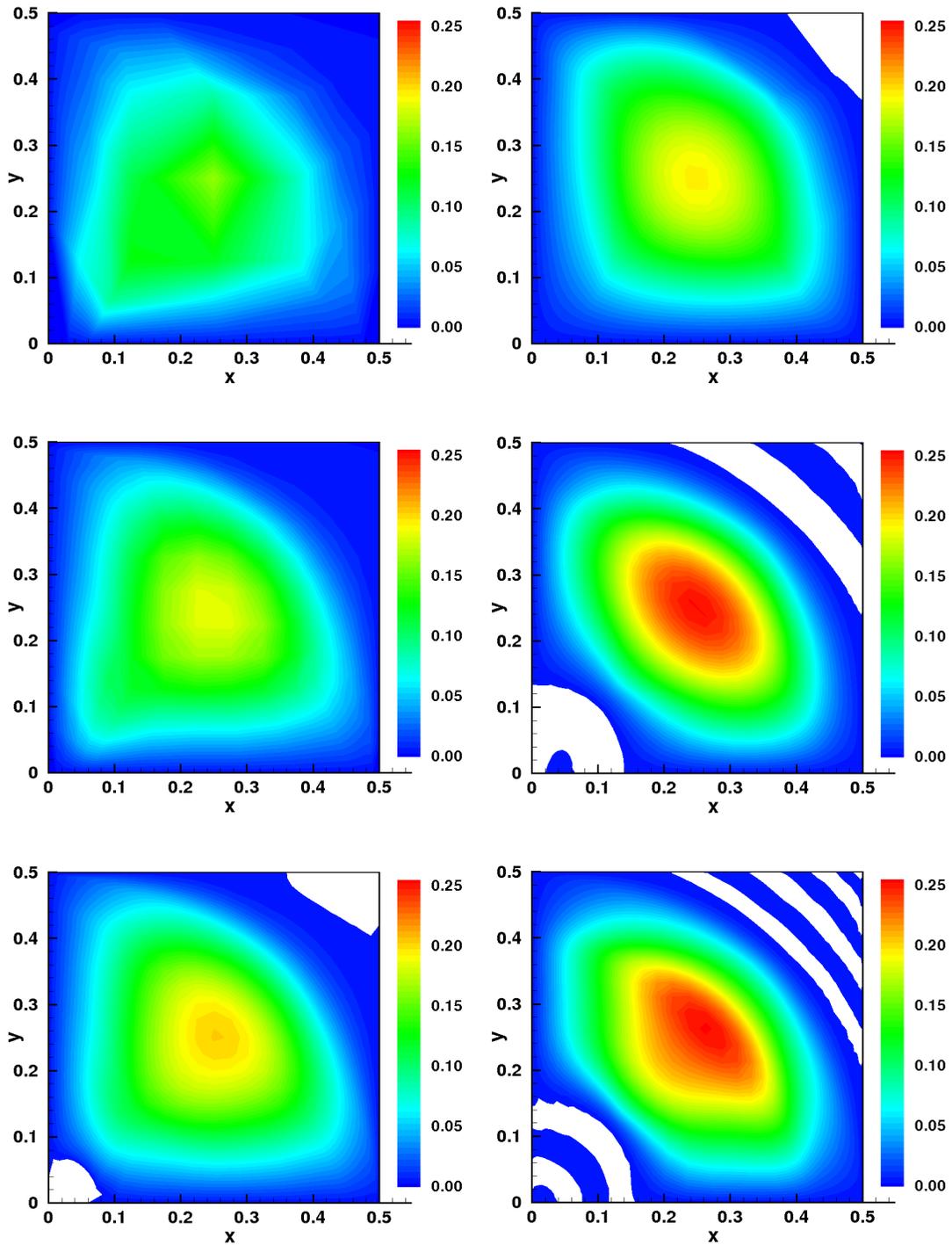

FIGURE 17. Non-uniform anisotropic media: This figure shows the concentration contours obtained using the *LS1 formulation* and *unstructured mesh* (which is shown in Figure 12). The left figures are using $h$-refinement, and the right figures are using $p$-refinement. The top figures are for $h/p = 2$, the middle figures are for $h/p = 5$, and the bottom two figures are for $h/p = 10$.

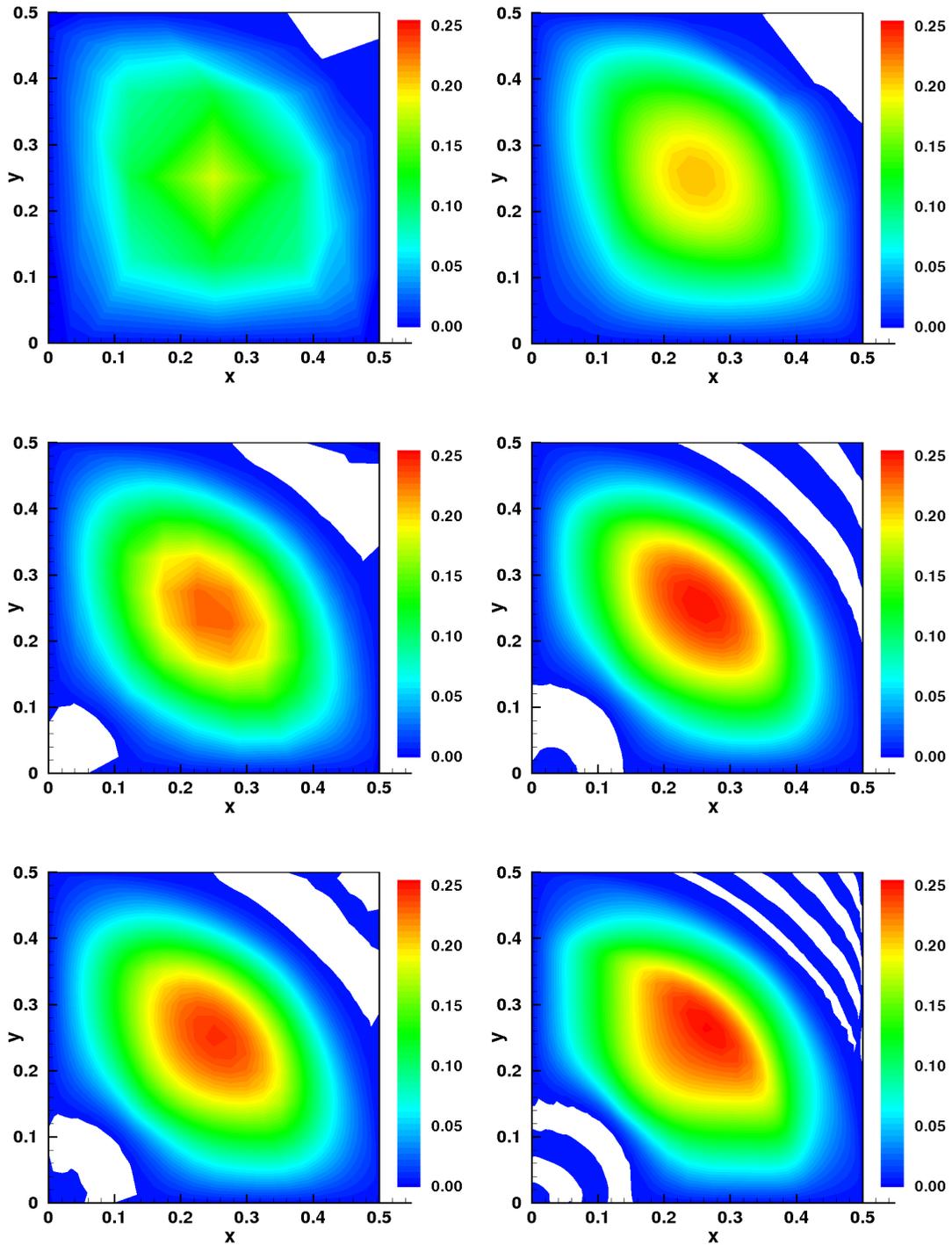

FIGURE 18. Non-uniform anisotropic media: This figure shows the concentration contours obtained using the *LS2 formulation* and *unstructured mesh* (which is shown in Figure 12). The left figures are using $h$-refinement, and the right figures are using $p$-refinement. The top figures are for $h/p = 2$, the middle figures are for $h/p = 5$, and the bottom two figures are for $h/p = 10$.

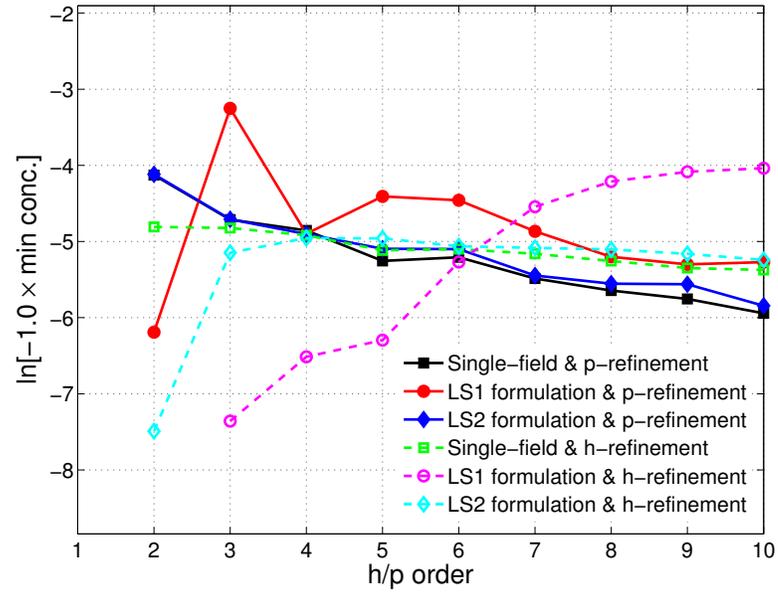

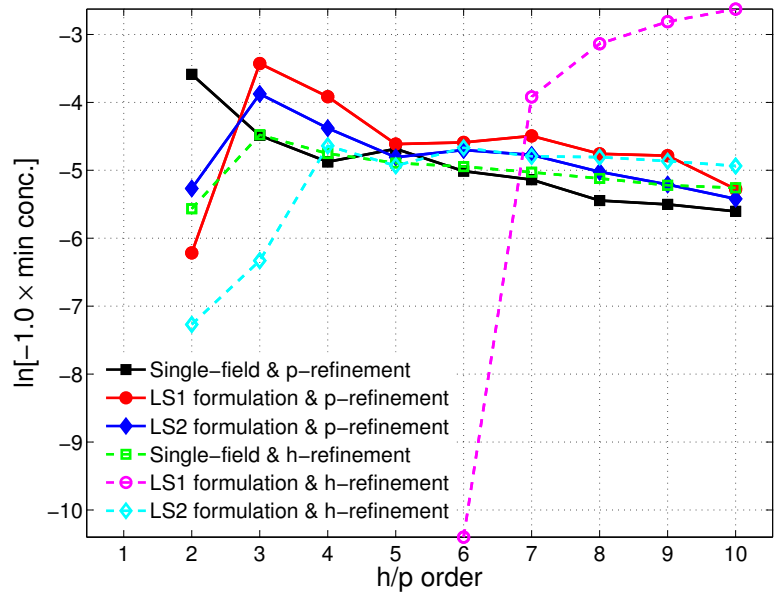

Figure 19. Non-uniform anisotropic media: This figure shows the minimum concentration for various order of polynomial approximation for single-field, LS1 and LS2 formulations using structured (top figure) and unstructured (bottom figure) meshes. According to the maximum principle, the concentration should be non-negative and should occur on the boundary. However, for all the three formulations the obtained minimum concentration is negative. The violation of the non-negative constraint did not vanish under either $p$-39or $h$-refinement. For better comparison, we have taken the $y$-axis to be the logarithm of the negative of the minimum concentration.

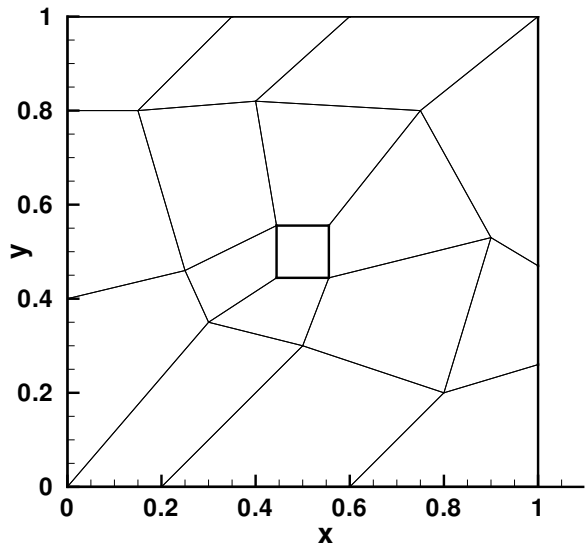

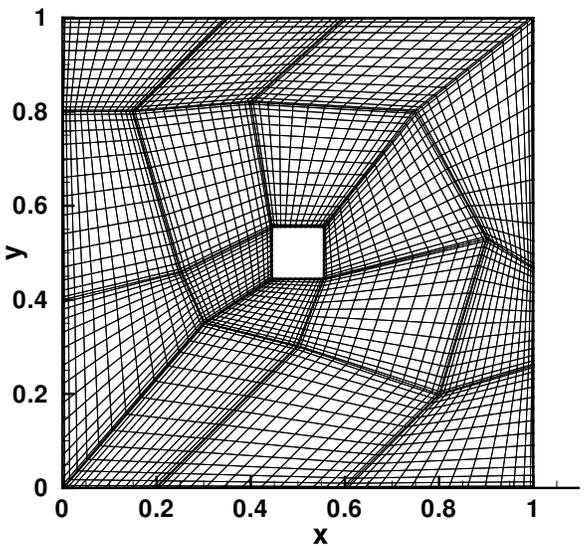

FIGURE 20. Anisotropic diffusion in a square domain with a hole: The top figure shows the mesh used in the numerical simulation, and the bottom figure shows the corresponding mesh used for visualization under high-order approximations. For more details on visualizing results using high-order approximations see subsection 5.1.

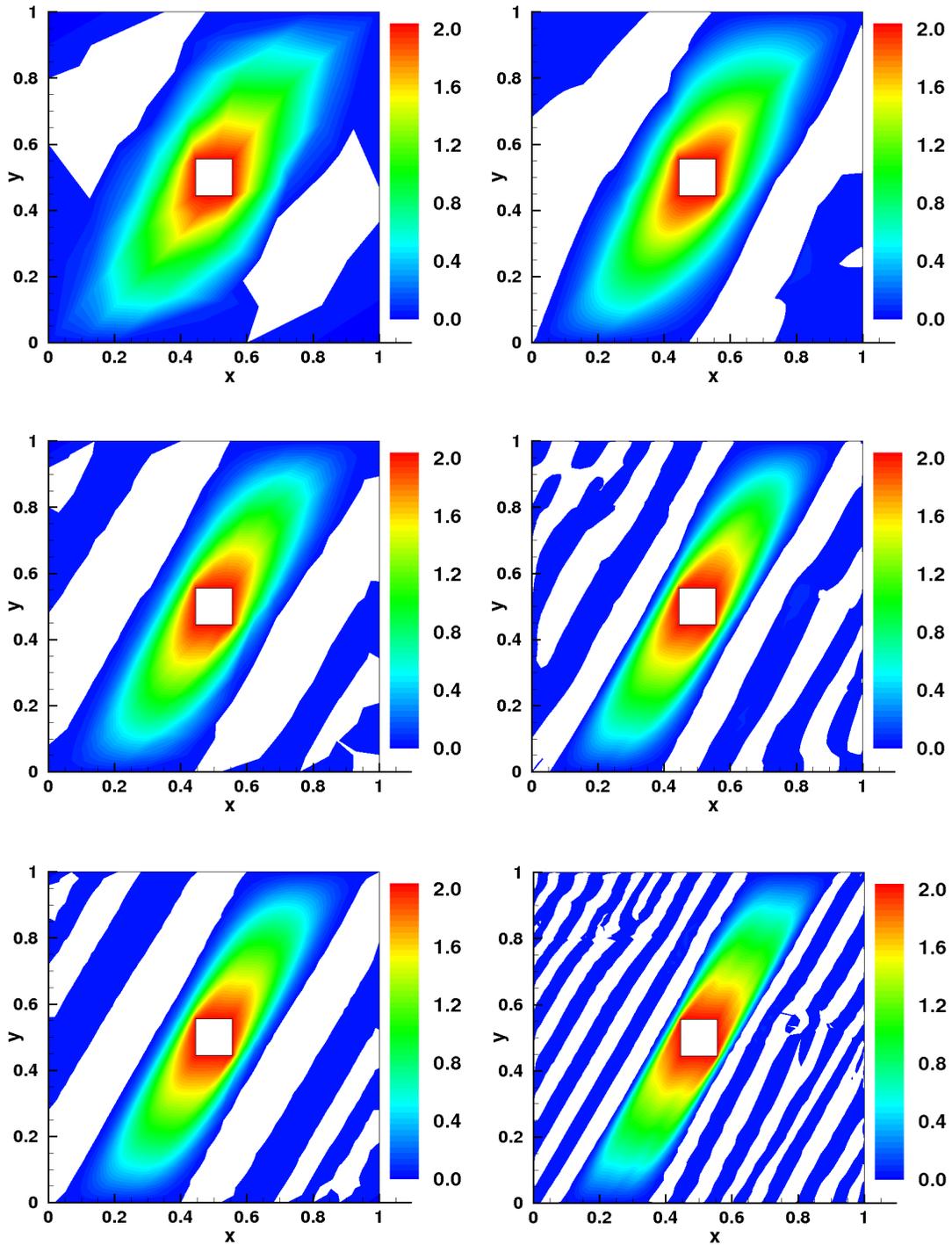

FIGURE 21. Anisotropic diffusion in a square domain with a hole: This figure shows the concentration contours obtained using the *single-field formulation* for $k_1 = 1$ and $k_2 = 10000$. The left figures are for $h$-refinement, and the right figures are $p$-refinement. The top figures are for $h/p = 2$, the middle figures are for $h/p = 5$, and the bottom two figures are for $h/p = 10$.

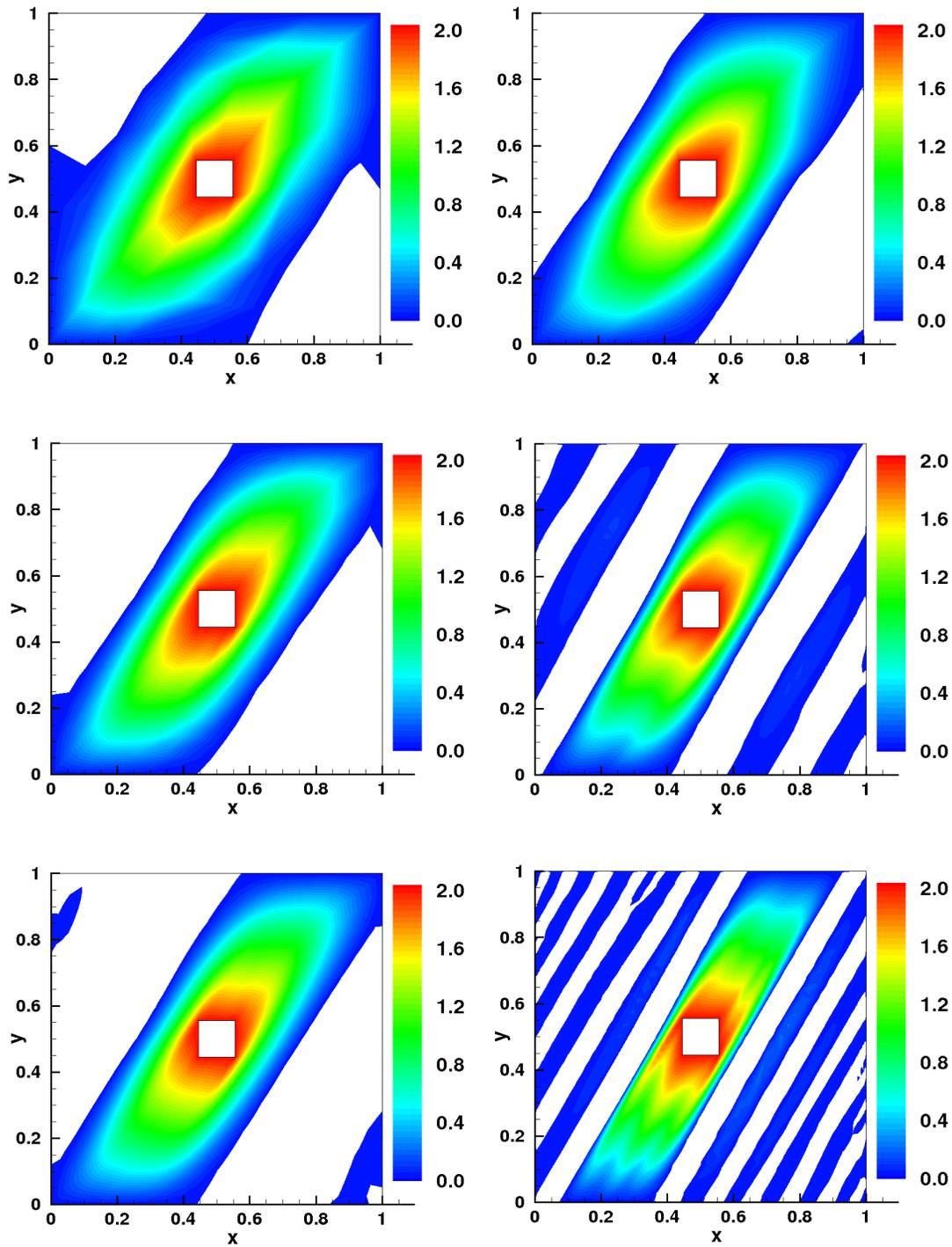

FIGURE 22. Anisotropic diffusion in a square domain with a hole: This figure shows the concentration contours obtained using the *LS1 formulation* for $k_1 = 1$ and $k_2 = 10000$. The left figures are using $h$-refinement, and the right figures using $p$-refinement. The top figures are for $h/p = 2$, the middle figures are for $h/p = 5$, and the bottom two figures are for $h/p = 10$.

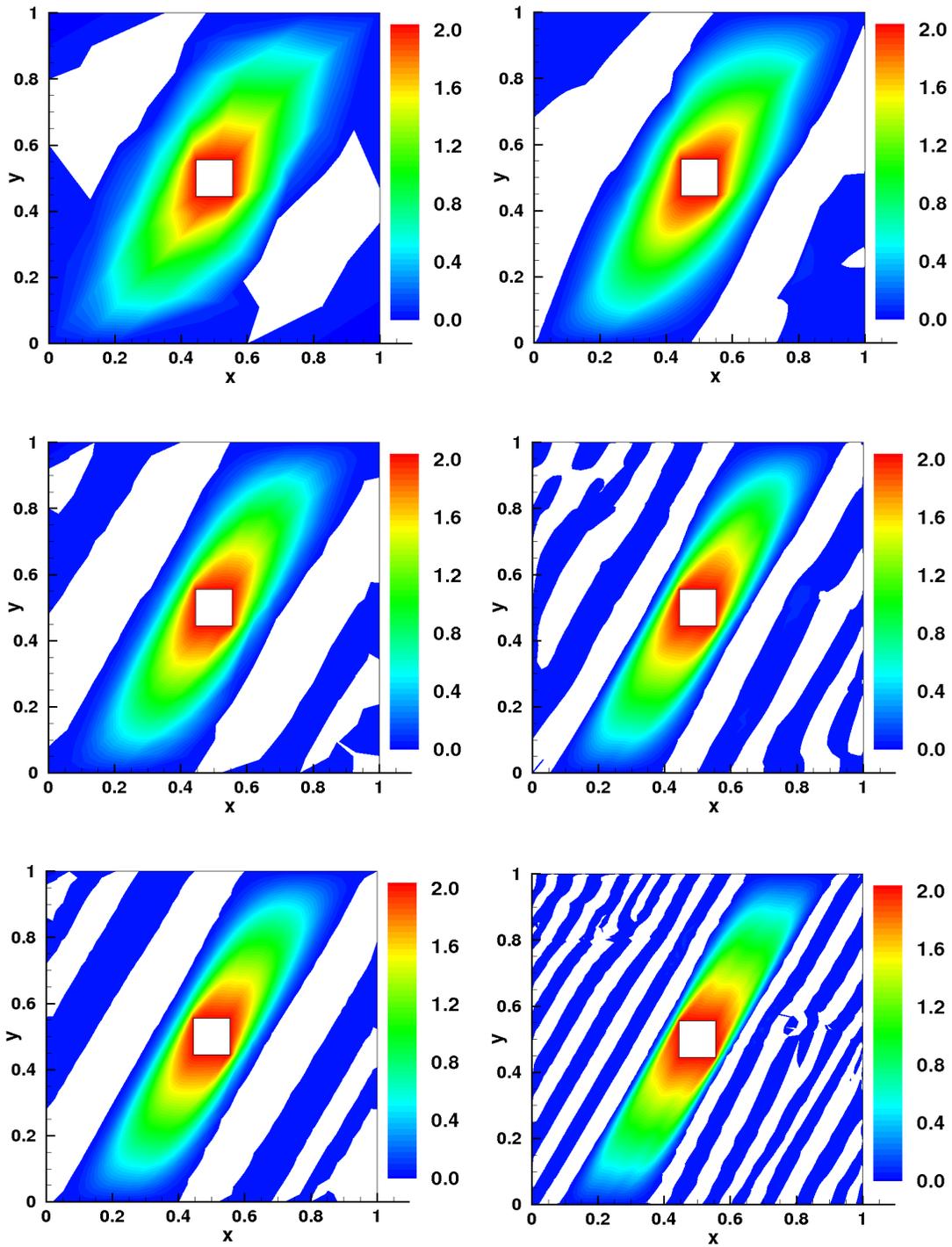

FIGURE 23. Anisotropic diffusion in a square domain with a hole: This figure shows the concentration contours obtained using the *LS2 formulation* for $k_1 = 1$ and $k_2 = 10000$. The left figures are for $h$-refinement, and the right figures are $p$-refinement. The top figures are for $h/p = 2$, the middle figures are for $h/p = 5$, and the bottom two figures are for $h/p = 10$.[43]

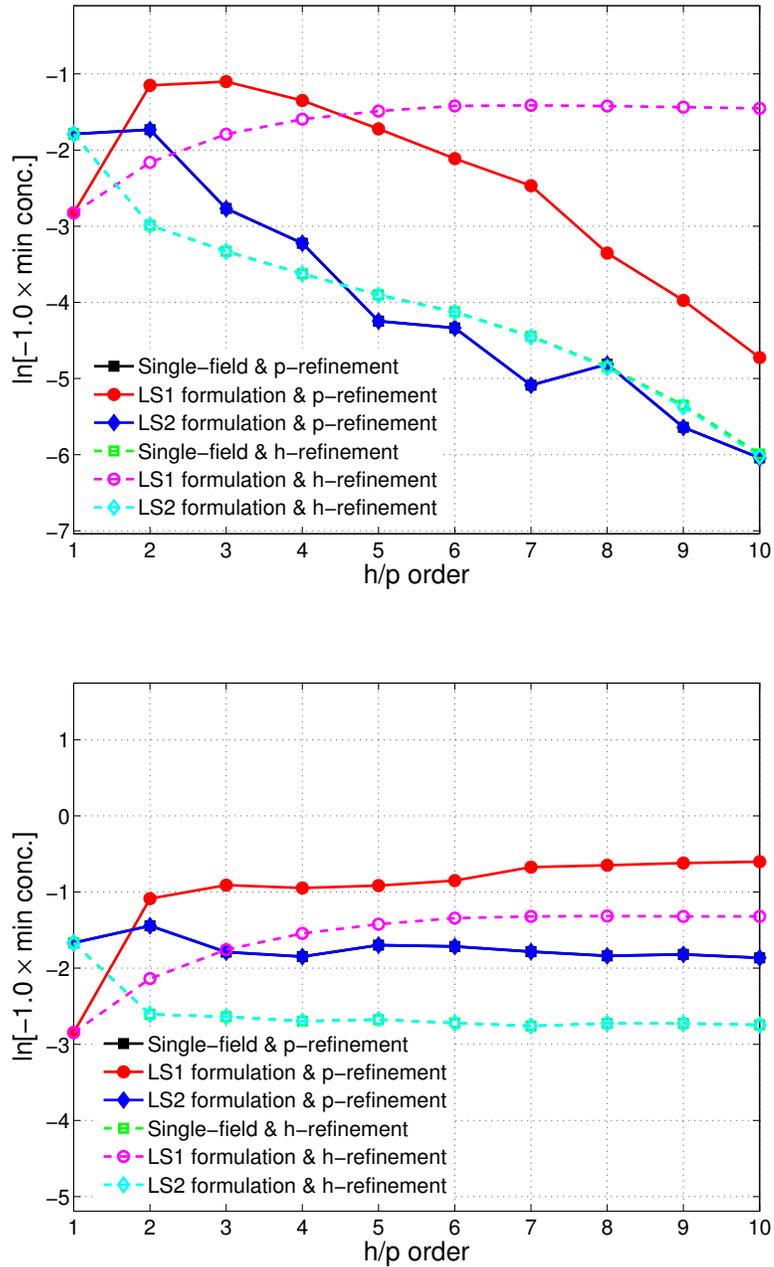

FIGURE 24. Anisotropic diffusion in a square domain with a hole: This figure shows the minimum concentration under single-field, LS1 and LS2 formulations for two different cases: $k_1 = 1$ and $k_2 = 100$ (top figure), and $k_1 = 1$ and $k_2 = 10000$ (bottom figure). We have shown the minimum concentration under various levels of $p$- and $h$-refinements. For all the three formulations and under various levels of $p$- and $h$-refinements, the minimum concentration is negative. For better comparison, we have taken the $y$-axis to be the logarithm of the negative of the minimum concentration. (Note that the minimum concentrations obtain using the single-field and LS2 formulations are approximately the same under both $p$- and $h$-refinements.)